\documentclass[letterpaper, 10 pt, journal]{IEEEtran} 	

\IEEEoverridecommandlockouts    % This command is only needed if 
\usepackage{amsthm}
\theoremstyle{definition}
\newtheorem{theorem}{Theorem}  %[section]
\newtheorem{lemma}[theorem]{Lemma}

\newtheorem{proposition}[theorem]{Proposition}
\newtheorem{problem}{Problem}

\newtheorem{assumption}{Assumption}
\newtheorem{example}{Example}

\newcommand{\mc}{\mathcal}

\newcommand{\trace}{\operatorname{Trace}}
\newcommand{\real}{\mathbb{R}}

\newcommand{\integnneg}{\mathbb{Z}_{\geq 0}}
\newcommand{\integpos}{\mathbb{Z}_{> 0}}

\newcommand{\naturalnneg}{\mathbb{N}_{\geq 0}}
\newcommand{\realpos}{\mathbb{R}_{> 0}}
\newcommand{\realnneg}{\mathbb{R}_{\geq 0}}

\newcommand{\inv}{{\negat 1}} 
\newcommand{\negat}{\scalebox{0.75}[.9]{\( - \)}}
\newcommand*{\QEDB}{\hfill\ensuremath{\square}}%  empty square

\newcommand*{\QEDBL}{\hfill\ensuremath{\blacksquare}}% Black square
\newcommand{\diag}{\operatorname{diag}}

\usepackage{pifont} % Needed for \ding 
\renewcommand{\ast}{\text{\ding{86}}}

\def\symmetric{\mathbb{S}}
\newcommand{\ftwo}{\tilde f}
\newcommand{\fone}{f}
\newcommand{\Lone}{L} % Lipschitz constant of \fone
\newcommand{\Ltwo}{\tilde L} % Lipschitz constant of \ftwo
\newcommand{\Lphi}{L_\phi} % Lipschitz constant of \phi
\newcommand{\fp}{z} % f perturbed in the algorithm

\newcommand{\sbs}[2]{{#1}_{\textup{#2}}}

\newcommand{\norm}[1]{\Vert #1 \Vert}

\DeclareMathAlphabet{\mymathbb}{U}{BOONDOX-ds}{m}{n}
\newcommand{\one}{\mathds{1}}
\newcommand{\zero}{\mymathbb{0}}

\newcommand{\reljust}[2]{\overset{\text{#2}}{#1}} % justify relationship
\newcommand*\mathcircled[1]{\tikz[baseline=(char.base)]{
    \node[shape=circle,draw,inner sep=1pt] (char) {\ensuremath{#1}};}}

\usepackage{enumerate}
\usepackage{graphicx,color}
\definecolor{lightblue}{RGB}{90,170,255}
\definecolor{myblue}{RGB}{0,90,160}
\usepackage[
    colorlinks=true,
    linkcolor=myblue,
    citecolor=myblue,
    urlcolor=myblue
]{hyperref}
\graphicspath{{./figures/}}
\usepackage{amsmath}
\usepackage{amssymb}
\usepackage{subfigure}
\usepackage{url}
\usepackage{cite}
\usepackage[vlined,ruled]{algorithm2e}
\usepackage{dsfont}
\usepackage{parskip}
\usepackage[shortlabels]{enumitem}
\usepackage{tikz}
\usepackage[nameinlink]{cleveref}

\title{\LARGE\textbf{Model-Free Aggregative Cooperative Optimization via 
Randomized Gradient-Free Minimization and Exploration Momentum
}}

\begin{document}

\author{Amir Mehrnoosh, Giuseppe Speciale, Riccardo Brumali, Giuseppe Notarstefano, Gianluca Bianchin
\thanks{%
The research is supported in part by F.R.S.-FNRS.
A.~Mehrnoosh and G.~Bianchin are with the ICTEAM institute at 
UCLouvain, Belgium. 
\texttt{\{\href{mailto:amir.mehrnoosh@uclouvain.be}{amir.mehrnoosh},\href{mailto:gianluca.bianchin@uclouvain.be}{gianluca.bianchin}\} @uclouvain.be}. 
R. Brumali and G. Notarstefano are with the Dept. of Electrical, Electronic and Information Engineering, Universit\`a di Bologna, Bologna, Italy. 
\texttt{\{\href{mailto:riccardo.brumali@unibo.it}{riccardo.brumali},\href{mailto:giuseppe.notarstefano@unibo.it}{giuseppe.notarstefano}\}@unibo.it}. G. Speciale was affiliated with Universit\`a di Bologna, Bologna, Italy, and the bulk of this work was conducted while he was a visiting student at UCLouvain.
} \hspace{-1cm}}

\maketitle

\thispagestyle{plain}\pagestyle{plain}

\begin{abstract}
Aggregative cooperative optimization problems arise in distributed 
decision-making settings where each agent's objective depends on its own 
decision 
as well as on an aggregate variable capturing global system behavior. Motivated 
by practical scenarios where gradient information is unavailable, this paper 
introduces a randomized gradient-free algorithm, named \texttt{ARGFree}, for 
solving such problems. 
\texttt{ARGFree} combines finite-difference gradient approximations with a set 
of tracking variables, emulating the behavior of a gradient-based method. We 
prove that \texttt{ARGFree} converges in expectation to an 
approximate optimizer, with the approximation error stemming from the use of a 
randomized gradient estimator. 
To enhance performance in high-dimensional settings, we further propose an improved variant, \texttt{ARGFree-EM}, which incorporates momentum in the exploration signals to smooth sudden fluctuations in the gradient exploration signals and thereby improve the accuracy of the underlying distributed tracking mechanism.
To the best of our knowledge, the class of \texttt{ARGFree} methods is the 
first in the literature capable of solving aggregating cooperative 
optimization problems without gradient information.

\end{abstract}

%%%%%%%%%%%%%%%%%%%%%%%%%%%%%%%%%%%%%%%%%%%%%%%%%%%%%%%%%%%%%%%%%%%%%%
%%%%%%%%%%%%%%%%%%%%%%%%%%%%%%%%%%%%%%%%%%%%%%%%%%%%%%%%%%%%%%%%%%%%%%%%%%%%%%%
\section{Introduction}
\label{sec:introduction}

In recent years, the widespread use of multi-agent systems has sparked 
increasing interest in solving optimization problems through distributed 
approaches. Decentralized decision-making schemes are particularly well-suited 
for applications involving large numbers of agents, or cases where information 
needs to maintained private. 
Representative examples include parameter estimation and detection, source 
localization in sensor networks, utility maximization, resource allocation, and 
multi-robot coordination. See the recent 
surveys~\cite{AN-JL:18,TY-XY-JW-etal:19,GN-IN-AC:19} for comprehensive overviews 
of the field. 

\begin{figure}[t]
\centering
\includegraphics[width=.85\columnwidth]{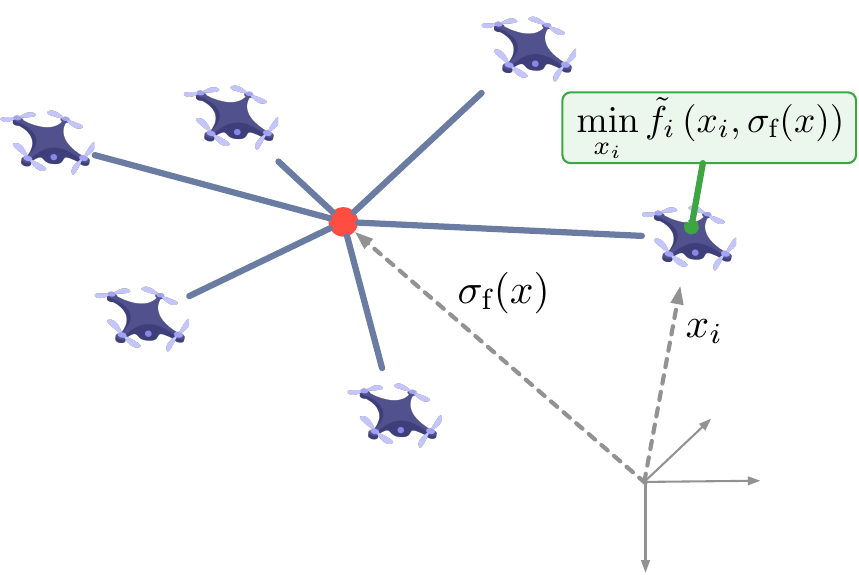}
\caption{
Illustrative application of the aggregative optimization framework 
developed in this work to robotic formation control problems.
Each robot $i$ controls its position $x_i,$ and seeks to optimize a 
local objective $\ftwo_i(x_i, \sbs{\sigma}{f}(x))$,  which captures 
both local positioning goals and coupling with the collective motion 
of the group, represented by the formation barycenter 
\( \sbs{\sigma}{f}(x) := \frac1N \sum_{i=1}^N x_i \). See 
Section~\ref{sec:illustrative_applications} for additional 
illustrative applications.
}
\vspace{-.2cm}
\label{fig:applications_robotics}
\end{figure}

A large body of work on cooperative distributed optimization focuses on the 
so-called \emph{consensus optimization} or \emph{federated learning} 
framework~\cite{AN-JL:18}. In this setting, agents aim to jointly solve 
optimization problems of the form
% $$\min_{x_1, \dots, x_N} \frac1N \sum_{i=1}^N \fone_i (x_i),$$ subject to 
% \( x_i = x_j \) for all \( i \neq j \), 
\begin{align*}
\min_{x_1, \dots, x_N} \frac1N \sum_{i=1}^N \fone_i (x_i), && \text{subject to:} \quad x_i = x_j, i \neq j, 
\end{align*}
and potentially additional problem-specific constraints.
In this formulation, \( \fone_i \) denotes the local loss function that agent 
\(i\) seeks to minimize, and the consensus constraints \( x_i = x_j \) ensures 
that all decision variables agree at convergence.
One key feature of this problem is that each $\fone_i$ depends only on the 
the local decision variable \( x_i \).
However, in many practical problems (such as feedback 
optimization~\cite{mehrnoosh2025optimization,GC-NM-GN:24}) local objective functions may 
depend not only on the agent’s own decision variable, but also on those of 
other agents. For example, in robotic formation problems, each agent may be 
interested in reaching configurations that depend not only on local 
objectives, but also on the barycenter of the group (see 
Fig.~\ref{fig:applications_robotics}). 
Moreover, the decision variables are not required to coincide at convergence, 
as in \(k\)-agreement problems~\cite{GB-MV-JC-ED:24-tac-b}.

These objectives have inspired the framework of \textit{aggregative 
cooperative optimization}, recently proposed in~\cite{XL-LX-YH:22}, 
where agents cooperatively solve optimization problems of the form
\[
\min_{x_1, \dots, x_N} \frac{1}{N} \sum_{i=1}^N 
\ftwo_i\!\left(x_i, \sbs{\sigma}{f}(x)\right),
\]
with \(\sbs{\sigma}{f}(x)\) denoting an aggregation function that depends on 
all decision variables \(x_1, \dots, x_N\).
Unfortunately, existing techniques to solve aggregative cooperative 
optimization problems assume that the agents have direct access to gradient 
(or subgradient) information about the local objective. In various 
applications, the relationship between the decision variables and the 
cost functions may be unknown, gradient information may be inaccessible, or 
the functions may not even be differentiable. 
Moreover, the dependence of $f$ on $\sbs{\sigma}{f}$ and the specific structure 
of $\sbs{\sigma}{f}$ further complicate the computation of gradients.
Motivated by this gap, in this work 
we introduce {\texttt{ARGFree}}: a distributed, gradient-free 
method that utilizes randomized (zeroth-order) finite-difference 
approximations to estimate gradients and solve aggregative cooperative 
optimization problems without gradients.

The \texttt{A}ggregative \texttt{R}andom \texttt{G}radient-\texttt{Free} 
(\texttt{ARGFree}) algorithm proposed here builds upon two main components:  
(i) a descent step, based on a forward-difference approximation of the 
gradient, which drives the optimization toward a minimizer of the aggregate 
loss function;  and (ii) a group of tracking variables, designed to estimate 
finite-difference approximations of the gradient. 
We demonstrate that the proposed method is capable of computing an approximate 
optimizer, where the approximation error stems from 
using a randomized gradient estimate in place of the exact gradient. 
Interestingly, our analysis reveals a trade-off between the rate of variation 
of the random exploration signal and the convergence rate of the tracking 
variables. To mitigate this effect, we propose a modified version of 
{\texttt{ARGFree}} with {\texttt{E}}xploration {\texttt{M}}omentum, called 
{\texttt{ARGFree-EM}}, which incorporates a damping term into the exploration 
signal to achieve improved accuracy.

\paragraph{Related work}
We classify the existing literature relevant to our work into two main 
categories: (i) the body of research on aggregative cooperative optimization, 
and (ii) gradient-free methods for distributed optimization.
\textit{Aggregative cooperative optimization:}
The aggregative cooperative optimization framework was 
introduced~\cite{XL-LX-YH:22}. In a departure from more  
classical {consensus optimization} problems~\cite{AN-JL:18}, this framework 
allows one to model a setting where the local decision need not to coincide at 
convergence, making it suitable for applications where coordination, rather 
than consensus, is the primary 
objective~\cite{mehrnoosh2025optimization,GC-NM-GN:24,GB-MV-JC-ED:24-tac-b}.
Our use of the wording \textit{``cooperative''} in combination with
\textit{``aggregative cooperative optimization''} is intended to distinguish our 
framework from that of \textit{aggregative games}~\cite{SL-PY-YH:17}, where 
each agent noncooperatively minimizes its own local objective while ignoring 
the effect of its decision on the rest of the network (and consequently, on the 
aggregate variable).
In a departure from aggregative games, this work is motivated by applications 
in feedback optimization~\cite{mehrnoosh2025optimization,GB-JC-JP-ED:21-tcns}, where the agents
act cooperatively to jointly optimize a common objective.
Online and constrained variants of the aggregative cooperative optimization 
problem have been studied in~\cite{XL-XY-LX:21}. Other notable contributions 
include~\cite{TW-PY:23}, which introduces a distributed Frank–Wolfe method, and 
the accelerated algorithms proposed in~\cite{LC-GW-XF-JZ-JC:24}. Particularly 
relevant to the present work are the recent works~\cite{GC-GN:22b,RB-GC-GN:25}, 
which harness learning-based techniques to handle uncertainty in the 
environment.
\textit{Gradient-free methods in distributed optimization:}
Although gradient-free techniques have a long history in optimization, their 
theoretical analysis was formalized only recently in~\cite{YN-VP:17}
in the centralized setting. In distributed settings, most existing works 
focus on the consensus optimization problem. Methods based on multi-point 
gradient estimators are studied in~\cite{YT-JZ-NL:20}, two-point estimators are 
proposed in~\cite{DH-MH-AG:19}, and single-point estimators are considered 
in~\cite{EM-M1:23}. Continuous-time algorithms have also been shown to be 
effective~\cite{JL-CT:12}. 
Further developments include analysis over time-varying graphs~\cite{YP-GH:19}, 
primal-dual approaches~\cite{DY-DH:14}, constrained stochastic 
problems~\cite{AS-SK:20}, and communication-imperfect 
settings~\cite{DW-ZJ-ZW-WW:17}. Accelerated variants are proposed 
in~\cite{OB-PM-EB:20}, while extremum-seeking based approaches are studied 
in~\cite{NM-GC-AT-GN:24}. Gradient-free methods tailored for games are 
presented in~\cite{YH-JH:24}.
To the best of the authors' knowledge, all existing gradient-free methods are 
designed for the framework of consensus optimization or aggregative games; 
solutions specifically tailored to the class of aggregative cooperative 
optimization problems remain unexplored.

\paragraph{Contributions} 

This paper features three main contributions. First, we 
introduce {\texttt{ARGFree}} (Algorithm~\ref{alg:ARGFree}), a randomized, 
gradient-free method designed to solve aggregative cooperative optimization 
problems. Unlike existing approaches for this class of problems, our method 
does not require gradient information; instead, it relies only on local 
function evaluations (as detailed in 
requirements~\ref{requirement_1}–\ref{requirement_2} in 
Section~\ref{sec:formulation}).
Second, we establish convergence bounds for the proposed algorithm 
(Theorem~\ref{thm:algorithm_convergence}), showing that its iterates converge 
in expectation to an approximate optimizer; the approximation error arises 
from the use of a randomized gradient estimator in place of the exact gradient.
Importantly, the radius of the convergence neighborhood can be controlled 
by appropriately tuning the algorithm’s parameters. 
Third, we propose {\texttt{ARGFree-EM}}
(Algorithm~\ref{alg:ARGFree_filtration}), which incorporates momentum in the 
random exploration signals to reduce errors in the tracking variables caused 
by rapidly-varying perturbation signals. This modification yields improved 
convergence bounds (Theorem~\ref{thm:algorithm_convergence_filtration}). 
Finally, we validate our theoretical findings through numerical simulations 
on robotic formation control problems (Section~\ref{sec:simulations}).

% \textit{(Organization)}
\paragraph{Organization}
Section~\ref{sec:formulation} formulates the problem and presents some 
technical preliminaries. Section~\ref{sec:algorithm_design} presents 
the {\texttt{ARGFree}} algorithm and establishes its convergence 
properties, while Section~\ref{sec:single_realization} 
introduces {\texttt{ARGFree-EM}}, and presents its convergence properties.
Section~\ref{sec:analysis} presents the analysis of the methods, 
Section~\ref{sec:illustrative_applications} discusses some 
illustrative applications of the framework, 
Section~\ref{sec:simulations} reports numerical results, and 
Section~\ref{sec:conclusion} concludes the paper. Finally, some proofs 
are reported in the appendix.

% \textit{(Notation)}
\paragraph{Notation}
We let $\integpos$ and $ \integnneg$ be the set of positive and non-negative 
integers, respectively; $\mathbb{R}^n,$ $n \in \integpos$, is the set of 
$n$-dimensional real vectors; $\realpos^n$ and $\realnneg$ denote the set of 
vectors in $\real^n$ with, respectively, positive and nonnegative coordinates.
$\symmetric^n$ is the space of $n\times n$ symmetric real matrices. 
We let $[n]:=\{1,2, \ldots, n\}$; we denote by 
$\operatorname{col}\left(v_1, \ldots, v_n\right)$ the column 
vector obtained by stacking the vectors $v_1, \ldots, v_n$. 
When $v_1, \ldots, v_n$ have the same dimension, 
$\bar v := \frac{1}{N} \sum_{i=1}^n v_i$ denotes their average.
The notation $\|\cdot\|, x^{\top}$, and $\langle x, y\rangle$ indicates the 
standard Euclidean norm (or induced matrix norm), the transpose of 
$x \in \mathbb{R}^n$, and standard inner product of $x, y \in \mathbb{R}^n$. 
The symbols $\one_n$ and $\zero_n$ denote the $n$-dimensional column vector 
of all ones and all zeros, respectively; $I_n$ is the identity matrix 
of dimension $n.$ Dimensions will be omitted when they are clear from the 
context.  We let $\rho(M)$ denote the spectral radius of a square matrix 
$M.$ The symbol $\otimes$ is the Kronecker product. \(E_u[\cdot]\) denotes the expectation with respect to the random variable \(u\).

%%%%%%%%%%%%%%%%%%%%%%%%%%%%%%%%%%%%%%%%%%%%%%%%%%%%%%%%%%%%%%%%%%%%%%
%%%%%%%%%%%%%%%%%%%%%%%%%%%%%%%%%%%%%%%%%%%%%%%%%%%%%%%%%%%%%%%%%%%%%%%%%%%%%%%
\section{Problem setting}
\label{sec:formulation}

In this section, we introduce the problem of interest and illustrate its 
relevance through representative applications.

\subsection{Problem formulation}
Consider a group of $N$ agents, where each agent $i \in [N]$ is associated 
with a local decision variable $x_i \in \mathbb{R}^{n_i}$ and a local loss 
function $\ftwo_i : \mathbb{R}^{n_i} \times \mathbb{R}^d \to \mathbb{R}$.
Let $x = \operatorname{col}(x_1, \ldots, x_N)$ and define 
$n := \sum_{i=1}^N n_i$.
We assume that each loss function $\ftwo_i$ depends not only on the local 
decision $x_i$ but also on a global quantity  
$\sbs{\sigma}{f}(x) \in \mathbb{R}^d$, referred to as the 
\emph{aggregative variable}, which aggregates information from all agents as
\begin{align}
\label{eq:sigma_x}
\sbs{\sigma}{f}(x)  :=\frac{1}{N} \sum_{i=1}^N \phi_i\left(x_i\right),
\end{align}
where $\phi_i: \mathbb{R}^{n_i} \rightarrow \mathbb{R}^d$. 
The agents aim to collaboratively solve the following optimization\footnote{
For notational clarity, throughout this work, the tilde notation (e.g., \(\ftwo\)) is used to denote functions that explicitly depend on the 
aggregative variable (i.e., functions of two arguments), whereas functions without a tilde (e.g., $\fone$)  represent the corresponding composite (single-argument) functions. 
The two are related through~\eqref{eq:rel_f_fhat}, as detailed shortly below.}  problem:
\begin{align}\label{eq:optimization_1}
\min _{x \in \mathbb{R}^n} \fone(x) & 
:= \frac1N \sum_{i=1}^N \ftwo_i\left(x_i, \sbs{\sigma}{f}(x)\right).
\end{align}
Equations~\eqref{eq:sigma_x}--\eqref{eq:optimization_1} formalize an 
\emph{aggregative cooperative optimization problem}~\cite{XL-LX-YH:22}, 
in which the agents aim to minimize a global objective defined as the average 
of their local costs. 
It is worth stressing that, unlike aggregative games~\cite{SL-PY-YH:17}, where 
agents optimize local costs noncooperatively (that is, neglecting the 
dependence of $\sbs{\sigma}{f}(x)$ on $x_i$), 
problem~\eqref{eq:optimization_1} is inherently cooperative as the dependence 
of $\sbs{\sigma}{f}(x)$ on $x_i$ is explicitly taken into account in the 
formulation~\eqref{eq:optimization_1}.

In this work, we focus on scenarios in which the agents operate under the 
following requirements (R):
\begin{enumerate}[label={(R\arabic*)},leftmargin=*]
\item \label{requirement_1}
The functions $\ftwo_i(\cdot, \cdot)$ and $\phi_i(\cdot)$, along with the 
decision variable $x_i$, need to be kept private to agent $i$ and, therefore, 
are not known by any other agent $j \neq i$.
\item \label{requirement_2}
Agent \( i \) does not have access to the analytic expressions 
(nor the derivatives) of \( \ftwo_i(\cdot, \cdot) \) or \( \phi_i(\cdot) \); 
instead, these functions can only be evaluated through oracle queries:
\begin{align*}
(x_i, \sigma) \mapsto \ftwo_i (x_i, \sigma), &&
x_i \mapsto \phi_i(x_i),
\end{align*}
with $x_i \in \real^{n_i}, \sigma \in \real^d.$ 
\end{enumerate}

Requirements~\ref{requirement_1}--\ref{requirement_2} arise naturally 
in many practical applications. Illustrative examples of such 
applications are discussed in Section~\ref{sec:illustrative_applications}.

We consider scenarios in which agents aim to solve 
problem~\eqref{eq:optimization_1} in a distributed manner, using only local 
communication and cooperative coordination. To this end, we assume that the 
agents exchange information with their neighbors; we model the communication 
topology using a directed graph $\mathcal{G}=(\mathcal{V}, \mathcal{E})$, 
where the node set $\mathcal{V}=\{1, \ldots, N\}$ models the agents and the 
edge set $\mathcal{E} \subset \mathcal{V} \times \mathcal{V}$ describes the 
communication links. 
We impose the following assumption on \(\mathcal{G}\) (see 
Section~\ref{sec:notation} for the adopted graph-theoretic notation).

\begin{assumption}[\textbf{\textit{Properties of the communication graph}}]
\label{as:graph}
The digraph $\mathcal{G}$ is strongly connected. Moreover, $\mc G$ admits an 
adjacency matrix $A$ that is doubly stochastic.
\QEDB\end{assumption}

Assumption~\ref{as:graph} is standard in the design of coordination schemes 
and distributed algorithms (e.g., consensus 
averaging)~\cite{VB-JH-AO-JT:05,RO-RM:04,WR-RWB-EMA:05}. 
Intuitively, the strong connectivity condition guarantees that information can 
(asymptotically) propagate from any node of \(\mathcal{G}\) to every other 
node, without necessarily implying that the graph is 
complete~\cite{VB-JH-AO-JT:05}. 
The doubly stochasticity requirement is, for instance, automatically satisfied 
when \(\mathcal{G}\) is additionally aperiodic or when each node possesses a 
self-loop~\cite{BG-JC:12}.
We also note that several procedures are available to construct matrices \(A\) 
that satisfy this assumption; see~\cite{VB-JH-AO-JT:05} for centralized 
methods and~\cite{BG-JC:12} for distributed ones. 
In the remainder, we let \(A \in \mathbb{R}^{N \times N}\) be a matrix as in 
Assumption~\ref{as:graph}.

We now formally state the problem studied in this work.

\begin{problem}[\textbf{\textit{Objective of this work}}]
\label{prob:main_problem}
Design a distributed algorithm, compatible with the graph topology $\mc G$, 
enabling the agents to cooperatively compute solutions to 
problem~\eqref{eq:optimization_1}, subject to 
requirements~\ref{requirement_1}--\ref{requirement_2}.
\QEDB\end{problem}

%%%%%%%%%%%%%%%%%%%%%%%%%%%%%%%%%%%%%%%%%%%%%%%%%%%%%%%%%%%%%%%%%%%%%%
%%%%%%%%%%%%%%%%%%%%%%%%%%%%%%%%%%%%%%%%%%%%%%%%%%%%%%%%%%%%%%%%%%%%%%%%%%%%%%%
\subsection{Preliminaries}
\label{sec:notation}
We now present basic properties used throughout the~paper.

% {\it (Basic notions on algebraic graph-theory)}
\paragraph{Notions on algebraic graph theory}
% The notation $\mathcal{G}=(\mathcal{V}, \mathcal{E})$ indicates a graph with 
% node set $\mathcal{V}=\{1, \ldots, N\}$ and edge set 
% $\mathcal{E} \subset \mathcal{V} \times \mathcal{V}$. 
For a digraph \(\mathcal{G} = (\mathcal{V}, \mathcal{E})\), we adopt the 
convention that an edge \((j, i) \in \mathcal{E}\) indicates that node 
$j$ is able to receive information from $i$ (or, equivalently, \(i\) transmits 
information to \(j\)). 
For a node $i \in \mc V$, we denote by 
$\mathcal{N}_i=\{j \in \mathcal{V}:(i, j) \in \mathcal{E}\}$ the set of  
agents that send information to $i$.
Matrix $A=\left[a_{i j}\right] \in \mathbb{R}^{N \times N}$ is said to be an 
{\it adjacency matrix} for $\mc G$ if it satisfies $a_{i j}>0$ if 
$(j, i) \in \mathcal{E}$, and $a_{i j}=0$ otherwise. 
$A$ is said to be {\it doubly stochastic} if $\sum_{j=1}^N a_{i j}=1$ 
and $\sum_{i=1}^N a_{i j}=1.$ 
We let $J:=\frac{1}{N} \one_N \one_N^{\top}\in \real^{N\times N}$ and  
$\mathcal{J}:=J \otimes I_d \in \real^{Nd \times Nd}$.
Notice that $\norm{\mc J - I}=1.$ 
With a slight abuse of notation, we denote by $\rho_A=\|A-J\|$ the operator 
norm of $A$ (which, in general, differs from its spectral radius \(\rho(A)\)).
Let $v = \operatorname{col}(v_1, \dots, v_N) \in \real^{Nd}$ with 
$v_i \in \real^d, i \in [N],$ and recall that 
$\bar v := \frac{1}{N} \sum_{i=1}^N v_i \in \real^d$ denotes the entries 
average. Recalling the well-known~\cite{RAH-CRJ:85}  property
$(A \otimes B)(C \otimes D)=A C \otimes B D$ for $A,B,C,D$ of 
suitable dimensions, we have 
\begin{align}\label{eq:kronone_calJ}
\one_N \otimes \bar v=\mc J v.
\end{align}
Given $A \in \real^{N \times N},$ we denote by 
$\mc A := A \otimes I_d \in \real^{Nd}.$
The following lemma is instrumental to our analysis.
\begin{lemma}[\!\!\cite{RAH-CRJ:85}]
Let $A \in \real^{N \times N}$ be a doubly stochastic matrix.
Then, $\rho_A$ satisfies $\rho_A<1$. Moreover, 
\begin{align}
\label{eq:inequalities_AJ}
\mathcal{A} \mathcal{J}=\mathcal{J} \mathcal{A}=\mathcal{J}, 
&& \|\mathcal{A} x-\mathcal{J} x\| \leq \rho_A\|x-\mathcal{J} x\|,
\end{align}
for any $x \in$ $\mathbb{R}^{N d}$.
\QEDB\end{lemma}

% Let $P \in \real^{k \times k}, k \in \integpos$ be a communication matrix
% associated with the graph $\mc G_P$. The matrix $P$ is said 
% \textit{irreducible} if its associated graph is strongly 
% connected. A matrix is said \textit{nonnegative} if all its entries are nonnegative.
% A (nonnegative) communication matrix is called row (or column) 
% \textit{stochastic} if all of its row-sums (or column-sums) are $1$. 
% An irreducible stochastic matrix $P$ is said \textit{primitive} if it has 
% only one eigenvalue with maximum modulus. 

% \begin{lemma}[\!\!\cite{RAH-CRJ:85}]
% Let $P \in \mathbb{R}^{n \times n}$ be an irreducible nonnegative matrix. 
% If $P$ has at least one nonzero diagonal entry, then $P$ is primitive.
% \QEDB\end{lemma}

% \begin{lemma}[\!\!\cite{RAH-CRJ:85}]
% Let $P \in \mathbb{R}^{n \times n}$ be an irreducible nonnegative matrix. 
% Then, the following properties hold: 1) $\rho(P)$ is an eigenvalue of $P$ 
% with algebraic multiplicity equal to $1$; 2) there exists a vector $x$ with 
% positive entries such that $P x=\rho(P) x.$
% \QEDB\end{lemma}

% \subsection{Gaussian smoothing}
% {\it (Preliminaries on Gaussian smoothing)}
\paragraph{Notions on Gaussian Smoothing-based  approximations}
Consider a function $\fone:\real^n \to \real$ and assume that, at each point 
$x \in \real^n$, it is differentiable along any direction. 
The \textit{Gaussian approximation} of $\fone(x)$ is defined as:
\begin{align}\label{eq:smoothed_f}
\fone_\delta(x)=\frac{1}{\kappa} \int_{\real^n} \fone(x+\delta u) \mathrm{e}^{-\frac{1}{2}\|u\|^2} \mathrm{~d} u,
\end{align}
where $\delta \in \realpos$ and 
$\kappa := \int_{\real^n} \mathrm{e}^{-\frac{1}{2}\|u\|^2} \mathrm{~d} u.$
Observe that, if we let $u \sim \mc N(0,\Sigma)$, then 
$\fone_\delta(x) \equiv E_u [\fone(x+\delta u)]$; in this case, 
$\kappa = (2 \pi)^{n / 2} (\operatorname{det} \Sigma)^{1 / 2}.$
Given $u \sim \mc N(0,\Sigma)$, we define the \textit{forward-difference 
gradient-free oracle} as follows:
\begin{align}\label{eq:grad_oracle}
g_{\delta}(x) = \frac{\fone(x+\delta u)-\fone(x)}{\delta} \Sigma^\inv u,
\end{align}
where $\delta \in \realpos$. 
Intuitively, in \eqref{eq:grad_oracle}, the vector $u$ can be interpreted as a 
random perturbation or directional probe that excites the function $\fone(x)$ to 
estimate its gradient. The following results are instrumental to our 
analysis.

\begin{lemma}[\!\!{\cite[Lem.~1]{YN-VP:17}}]
\label{lem:expectation_norm_bound}
Let $p \in \integnneg$ and 
$M_p := \frac{1}{\kappa} \int_{\real^n}\|u\|^p
\mathrm{e}^{-\frac{1}{2}\|u\|^2} \mathrm{~d} u.$ Then, $M_0=1, M_1=n, $ and 
\begin{align}\label{eq:expectation_norm_bound}
M_p \leq \begin{cases}n^{p / 2}, & \text { if } p \in[0,2], \\ 
(p+n)^{p / 2}, & \text { if } p> 2.
\end{cases}
\end{align}
\QEDB\end{lemma}

\begin{lemma}[\!\!{\cite[Thm.~2]{YN-VP:17}}]
Suppose $\fone: \real^n \to \real$ is convex and Lipschitz continuous. Then, 
\begin{align}
E_u[g_\delta (x)] &= \nabla \fone_\delta (x) \label{eq:Nesterov21}
\end{align}
for any $x \in \real^n$
\QEDB\end{lemma}

\begin{subequations}
\begin{lemma}[\!\!\!{\cite[Thm.s~1 and 4, Lemmas~4 and 5]{YN-VP:17}}] \label{lem:gaussianSmooth}
Suppose $\fone: \real^n \to \real$ is convex and $\Lone_1$-smooth. 
Then, for any $x \in \real^n$, 
\begin{align}
\hspace{-1cm}
\vert \fone_\delta(x) - \fone(x) \vert &\leq \frac{\delta^2}{2}\Lone_1 n, \label{eq:Nesterov19}\\
\fone_\delta(x) - \fone(x)  &\geq 0, \label{eq:Nesterov11}\\
\Vert \nabla \fone(x) \Vert^2 &\leq 2\Vert \nabla \fone_\delta(x)  \Vert^2 + \frac{\delta^2}{2}\Lone_1^2(n+6)^3, \label{eq:Nesterov29}\\
E_u[ \| g_\delta (x) \|^2 ] &\leq \frac{\delta^2}{2} (n+6)^3 \Lone_1^2+2(n+4)\|\nabla \fone(x)\|^2,  \label{eq:Nesterov35}\\
E_u[\Vert  g_\delta (x) \Vert^2] &\leq 4(n+4) \Vert \nabla \fone_\delta (x) \Vert^2 + 3\delta^2\Lone_1^2(n+4)^3. \label{eq:Nesterov37}
\end{align}
\QEDB\end{lemma}
\end{subequations}

%%%%%%%%%%%%%%%%%%%%%%%%%%%%%%%%%%%%%%%%%%%%%%%%%%%%%%%%%%%%%%%%%%%%%%
%%%%%%%%%%%%%%%%%%%%%%%%%%%%%%%%%%%%%%%%%%%%%%%%%%%%%%%%%%%%%%%%%%%%%%%%%%%%%%%
\section{The \texttt{ARGFree} algorithm: design and convergence guarantees}
\label{sec:algorithm_design}
In this section, we introduce an iterative algorithm for solving 
Problem~\ref{prob:main_problem} and establish its convergence guarantees.

\subsection{The \texttt{ARGFree} algorithm}
\label{sec:algorithm_description}
The proposed method, called \texttt{A}ggregative \texttt{R}andom 
\texttt{G}radient-\texttt{Free} (\texttt{ARGFree}) algorithm, is presented in 
Algorithm~\ref{alg:ARGFree}. 
The algorithm is structured as follows: each agent $i \in [N],$ 
updates its local decision variable $x_i^k$ using a forward-difference 
optimization scheme (see line 1), driven by the term 
$\alpha \frac{p_i^k - \fp_i^k}{\delta} \cdot u_i^k,$  where $\alpha>0$ is 
the algorithm's stepsize and $\delta>0$ is a tunable parameter (hereafter 
called \textit{smoothing ratio}---see~\eqref{eq:smoothed_f}).
In this operation, the vector $u_i^k \in \real^{n_i}$ describes a random 
perturbation direction, while the variables $p_i^k$ and $\fp_i^k$ model 
local estimates for, respectively, the quantities $\fone(x_k+\delta u_k)$ and 
$\fone(x_k).$ 
To update the estimates $p_i^k$ and $\fp_i^k$, the algorithm uses two 
dynamic consensus tracking schemes ({see lines 5 and 6}) driven, locally, by 
the signals 
$\ftwo_i(x_i^{k}+\delta u_i^{k}, \sbs{\sigma}{f}(x^{k}+\delta u^{k}))$ and
$\ftwo_i(x_i^{k}, \sbs{\sigma}{f}(x^{k}))$, respectively. 
Finally, because evaluating the quantities 
$\sbs{\sigma}{f}(x^{k}+\delta u^{k})$  and $\sbs{\sigma}{f}(x^{k})$ would 
require global knowledge (through knowledge of the global vectors $x^{k}$, 
$u^{k}$ as well as of the functions $\phi_j(\cdot)$), these quantities are 
replaced by the local estimates $s_i^k$ and $\sigma_i^k$ in lines 5--6. 
In other words, $s_i^k$ and $\sigma_i^k$ are local proxies for, 
respectively, 
$\sbs{\sigma}{f}(x^{k}+\delta u^{k})$  and $\sbs{\sigma}{f}(x^{k}),$ and are
estimated through a dynamic consensus tracking scheme ({see lines 3--4}).

\begin{algorithm}[t]
\DontPrintSemicolon
\caption{\texttt{ARGFree} (agent $i$)}
\label{alg:ARGFree}
\KwData{Parameters $\alpha, \delta \in \realpos,$}
\textbf{Initializations:} Set
$k=0,  u_{i}^{0} \sim \mc N(0, I_{n_i})$, $x_{i}^{0} \in \real^{n_i},$ 
\begin{align*}
\sigma_{i}^{0} &= \phi_i(x_{i}^{0}), &
s_{i}^{0} &= \phi_i(x_{i}^{0} + \delta u_{i}^{0}), \\
\fp_{i}^{0} &= \ftwo_i(x_{i}^{0}, \sigma_{i}^{0}), &
p_{i}^{0} &= \ftwo_i(x_{i}^{0} + \delta u_{i}^{0}, s_{i}^{0})
\end{align*}

%
%\nl At each time $k\in \naturalnneg:$\\
%
\textbf{Optimization variable update:}\;
%\nl Update $x_i^k$ as:
\nl 
\begin{algomathdisplay}
x_i^{k+1} = x_i^k - \alpha \frac{p_i^k - \fp_i^k}{\delta} u_i^k
\end{algomathdisplay}\;
\textbf{Tracking variables update:}\;
\nl Generate $u_{i}^{k+1} \sim \mc N(0, I_{n_i})$ and update:\;
\begin{algomathdisplay}
\begin{aligned}
\nl \sigma_i^{k+1} &= \sum_{j\in \mc N_i} a_{ij} \sigma_j^k + \phi_i(x_i^{k+1})
- \phi_i(x_i^{k})\;\\
\nl s_i^{k+1} &= \sum_{j\in \mc N_i} a_{ij} s_j^k + \phi_i(x_i^{k+1} + \delta u_i^{k+1}) - \phi_i(x_i^{k} + \delta u_i^{k})\;\\
\nl \fp_i^{k+1} &= \sum_{j\in \mc N_i} a_{ij} \fp_j^k + \ftwo_i(x_i^{k+1},\sigma_i^{k+1}) - \ftwo_i(x_i^{k}, \sigma_i^k)\;\\
\nl p_i^{k+1} &= \sum_{j\in \mc N_i} a_{ij} p_j^k + \ftwo_i(x_i^{k+1} + \delta u_i^{k+1}, s_i^{k+1}) \\
& \hspace{3cm} - \ftwo_i(x_i^{k}+\delta u_i^{k}, s_i^k)\;\\
\end{aligned}
\end{algomathdisplay}\;
\nl Transmit $\sigma_i^{k+1}$, $s_i^{k+1}$, $\fp_i^{k+1},$ and $p_i^{k+1}$ to neighbors\;
\nl Set $k \gets k+1$ and go to line 1\;
\KwResult{$x_i^k$, estimate for the optimizer of~\eqref{eq:optimization_1}}
\end{algorithm}

\subsection{Convergence guarantees for \texttt{ARGFree}}

In this section, we establish convergence guarantees for~\texttt{ARGFree}. 
Before proceeding, we introduce some notation that is instrumental to state 
our results. We define $\ftwo :\real^n \times \real^{Nd} \to \real$ as
\begin{align*}
\ftwo(x,\sigma) := \frac1N \sum_{i=1}^N \ftwo_i\left(x_i, \sigma_i \right),
\end{align*}
where, for a vector $\sigma \in \real^{Nd},$ we used the notation
$\sigma = \operatorname{col}(\sigma_1, \dots, \sigma_N), \sigma_i \in \real^d.$
Observe that, by~\eqref{eq:sigma_x}--\eqref{eq:optimization_1},
the following identity holds
\begin{align}\label{eq:rel_f_fhat}
\ftwo(x, \one \otimes \sbs{\sigma}{f}(x))=\fone(x).
\end{align}

Next, we  introduce the vector notation:

\scalebox{.95}{\parbox{\linewidth}{
\begin{align}\label{eq:x_vector_notation}
x_k&=\operatorname{col}\left(x_{1}^{k}, \ldots, x_{N}^{k}\right) \in \real^n,\nonumber &
u_k&=\operatorname{col}\left(u_{1}^{k}, \ldots, u_{N}^{k}\right) \in \real^n,\nonumber \\
\sigma_k &=\operatorname{col}\left(\sigma_{1}^{k}, \ldots, \sigma_{N}^{k}\right) \in \real^{Nd},\nonumber &
s_k&=\operatorname{col}\left(s_{1}^{k}, \ldots, s_{N}^{k}\right) \in \real^{Nd},\nonumber \\
\fp_k&=\operatorname{col}\left(\fp_{1}^{k}, \ldots, \fp_{N}^{k}\right) \in \real^N, &
p_k&=\operatorname{col}\left(p_{1}^{k}, \ldots, p_{N}^{k}\right) \in \real^N.
\end{align}}}

Using this notation, the updates of Algorithm~\ref{alg:ARGFree} can be 
written in vector form\footnote{Given two vectors, 
$v = \operatorname{col}(v_1, \dots, v_N) \in \real^N, v_i \in \real,$ and 
$u =\operatorname{col}(u_1, \dots , u_N) \in \real^n, u_i \in \real^{n_i},$ we 
denote by $v \odot u= (v_1 u_1, \dots , v_N u_N) \in \real^n$ their  
entrywise product.} as:

\scalebox{.97}{\parbox{\linewidth}{
\begin{subequations}\label{eq:alg_vector_form}
\begin{align}
x_{k+1} &= x_k - \frac{\alpha}{\delta} (p_k-\fp_k) \odot u_k,\\
\sigma_{k+1} &=\mathcal{A} \sigma_k + \sbs{\phi}{v}\left(x_{k+1}\right)-
\sbs{\phi}{v} \left(x_k\right),\label{eq:sigma}\\
s_{k+1} &=\mathcal{A} s_k+\sbs{\phi}{v}\left(x_{k+1}+\delta u_{k+1}\right)-\sbs{\phi}{v}\left(x_k+\delta u_k\right),\label{eq:s}\\
\fp_{k+1} &=\mathcal{A} \fp_k+ \sbs{\ftwo}{v} \left(x_{k+1}, \sigma_{k+1}\right)- \sbs{\ftwo}{v} \left(x_k, \sigma_k\right),\label{eq:phi}\\
p_{k+1} &=\mathcal{A} p_k+\sbs{\ftwo}{v}\left(x_{k+1}+\delta u_{k+1}, s_{k+1}\right)-\sbs{\ftwo}{v}\left(x_k+\delta u_k, s_k\right),\label{eq:p}
\end{align}
\end{subequations}}}

% \ammargin{We also can use $x_{k+1} = x_k - \frac{\alpha}{\delta} \left( (p_k - \fp_k) \otimes I_{\text{block}} \right) u_k$, $I_{\text{block}} = \mathrm{blkdiag}(I_{n_1}, \dots, I_{n_N})$ in eq. (16a). 
% \red{[GB: this does not seem correct: $I_{\text{block}}$ will be an 
% $n \times n$ matrix, and thus  $(p_k - \fp_k) \otimes I_{\text{block}}$ 
% will be $Nn\times n.$ It works perhaps if all $u_i^k$ have the same dimension (the $n_i$ are equal).]}}

where $\sbs{\phi}{v}:\real^n \to \real^{Nd}$ and
$\sbs{\ftwo}{v}:\real^n \times \real^{N d}\to \real^N$ are:
\begin{align}
\sbs{\phi}{v}(x) &:=\operatorname{col}\left( \phi_1 (x_1), \dots,   \phi_N(x_N) \right), \notag\\
\sbs{\ftwo}{v}(x,\sigma) &
:=\operatorname{col}( \ftwo_1(x_1,\sigma_1), \dots,  \ftwo_N(x_N,\sigma_N)).
\end{align}

We begin with the following result, which formalizes the tracking properties of 
the dynamic consensus tracking scheme used in~\eqref{eq:sigma}--\eqref{eq:p} 
(or, lines 3--6 of Algorithm~\ref{alg:ARGFree}). 
Recall that, for a vector \( v = \operatorname{col}(v_1, \dots, v_N)\), the 
notation $\bar{v} := \frac{1}{N} \sum_{i=1}^N v_i$ denotes the average of its 
elements (see Section~\ref{sec:notation}).

\begin{proposition}[\textbf{\textit{Properties of the tracking variables}}]
\label{prop:tracking_variables}
Suppose Assumptions~\ref{as:graph}–\ref{as:objective} hold. Then, the states 
\((\sigma_k, s_k, \fp_k, p_k)\) of~\eqref{eq:alg_vector_form} satisfy:
\begin{align}\label{eq:averages_identities}
\bar \sigma_k &= \sbs{\sigma}{f}(x_k), &
\bar s_k &= \sbs{\sigma}{f}(x_k+\delta u_k), \notag\\
\bar \fp_k &= \ftwo(x_k,\sigma_k), &
\bar p_k &= \ftwo(x_k+\delta u_k, \sigma_k),
\end{align}
at every  $k \in \integnneg.$ 
\QEDB\end{proposition}

The proof of this claim is presented in 
Appendix~\ref{sec:proof_lem_tracking_variables}.
Proposition~\ref{prop:tracking_variables} establishes that the averages 
(computed across the agents) of the state variables 
\(\sigma_k, s_k, \fp_k, p_k\) 
coincide, respectively, with the quantities 
\( \sbs{\sigma}{f}(x_k) \), 
\( \sbs{\sigma}{f}(x_k + \delta u_k) \), 
\( \ftwo(x_k, \sigma_k) \), and 
\( \ftwo(x_k + \delta u_k, \sigma_k) \). 
This result implies that the averages across the network of the algorithm’s 
state variables \(\sigma_k, s_k, \fp_k,\) \(p_k\) respectively track the 
quantities they are designed to represent (see the interpretation of the state 
variables in Section~\ref{sec:algorithm_description}).

Next, we state an instrumental lemma that will be used to establish the 
convergence properties of Algorithm~\ref{alg:ARGFree}.

\begin{lemma}[\textbf{\textit{Contraction of randomized descent under 
forward-difference gradient approximation}}]
\label{lem:bound_twopoint_alg}
Let $\fone: \real^n \to \real$ be an $\Lone_1$-smooth and $\mu$-strongly 
convex function, and let $x^*$ be such that $\nabla \fone(x^*)=0.$ 
Let $u \sim \mc N(0,\Sigma)$ and consider the forward-difference gradient-free 
oracle (cf.~\eqref{eq:grad_oracle}):
\begin{align*}
g_{\delta}(x) = \frac{\fone(x+\delta u)-\fone(x)}{\delta} \Sigma^\inv u.
\end{align*}
Then, for any $x \in \real^n,$ the following inequality holds:
\begin{align}
\label{eq:contraction_forward_difference}
E_{u} [\norm{x - \alpha g_\delta(x)- x^*}] \leq 
\sqrt{1-\beta_1^\alpha} \norm{x - x^*} + \beta_2^\alpha,
\end{align}
where
\begin{align}\label{eq:def_betas}
\beta_1^\alpha  &:= \alpha \mu (1 - 2 \alpha (n+4)\Lone_1),\\
\beta_2^\alpha &:= \sqrt{\alpha \delta^2 \Lone_1 (n + \frac{\alpha}{2} (n+6)^3\Lone_1)}.
\tag*{\QEDB}
\end{align}
\end{lemma}
The proof of this claim is provided in 
Appendix~\ref{sec:proof_lem_bound_twopoint_alg}. 
The bound~\eqref{eq:contraction_forward_difference} concerns iterative 
algorithms of the form \( x_{k+1} = x_k - \alpha g_\delta(x_k) \), which can be 
interpreted as a (centralized) gradient-descent-type method that employs the 
forward-difference approximation~\eqref{eq:grad_oracle} in place of the exact 
gradient. 
The estimate~\eqref{eq:contraction_forward_difference} establishes 
that these iterates are contractive with respect to the optimizer 
\( x^* \). 
Specifically, the contraction occurs with rate \( \sqrt{1 - \beta_1^\alpha} \), 
up to a neighborhood of radius \( \beta_2^\alpha \). 
Here, convergence to an inexact point arises from employing a 
randomized gradient estimate instead of the exact gradient.
The convergence rate and accuracy of the algorithm depend on several 
optimization parameters; in particular, note that 
\( \beta_1^\alpha < 1 \) requires a stepize 
\( \alpha < \tfrac{1}{2(n+4)\Lone_1} \), and that the radius of the 
convergence neighborhood can be controlled (i.e., made arbitrarily 
small) by reducing the smoothing ratio \( \delta \).
% Lemma~\ref{lem:bound_twopoint_alg} establishes that the iterates of 
% $x_{k+1} = x_k - \alpha g_\delta(x_k)$---interpreted as a 
% gradient-descent-type algorithm employing the forward-difference 
% approximation~\eqref{eq:grad_oracle} for the gradient--are 
% contractive with respect to the optimizer \( x^* \), 
% with a contraction factor \( 1 - \beta_1^\alpha \), 
% up to a neighborhood of radius \( \beta_2^\alpha \).
% Note that ensuring contraction requires \( \beta_1^\alpha > 0 \), which holds 
% when the stepsize satisfies \( \alpha < \frac{1}{2(n + 4)\Lone_1} \).

Motivated by the statement of Lemma~\ref{lem:bound_twopoint_alg}, we impose 
the following requirements on the optimization 
problem~\eqref{eq:optimization_1}.

\begin{assumption}[\textbf{\textit{Properties of the loss functions}}]
\label{as:objective}
For all \( i \in [N] \), the following statements hold:
\begin{enumerate}[label={(A\ref{as:objective}\textup{\alph*})},leftmargin=*]
\item \label{as:objective_a}
The function $\fone(x)$ is Lipschitz smooth and 
$\mu$-strongly convex. We denote by $\Lone_0, \Lone_1 >0$ constants such 
that 
$\norm{\fone(x)-\fone(x^\prime)}\leq \Lone_0 \norm{x-x^\prime}$  and 
$\norm{\nabla \fone(x)- \nabla \fone(x^\prime)}\leq \Lone_1 \norm{x-x^\prime}$, 
$\forall x, x^\prime \in \real^n.$

\item \label{as:objective_b}
The functions $\phi_i(x_i)$ are Lipschitz continuous. We let $\Lphi>0$ be such 
that 
$\norm{\phi_i(x_i)-\phi_i(x_i^\prime)}\leq \Lphi \norm{x_i-x_i^\prime},$
$\forall x_i, x_i^\prime \in \real^{n_i}.$

\item \label{as:objective_c}
Each function $\ftwo_i$ is Lipschitz continuous. We let $\Ltwo_{0,i}>0$ be 
such that $\norm{\ftwo_i(x_i,\sigma)-\ftwo_i(x_i^\prime,\sigma^\prime)}\leq \Ltwo_{0,i} \norm{\begin{bmatrix} x_i\\ \sigma \end{bmatrix} -
\begin{bmatrix} x_i^\prime\\ \sigma^\prime \end{bmatrix}},$ 
$\forall x_i, x_i^\prime \in \real^{n_i}, \sigma, \sigma^\prime \in \real^d.$ Moreover, we let $\Ltwo_0:= \max_i \Ltwo_{0,i}.$
\QEDB\end{enumerate}
\end{assumption}

Note that strong convexity is required only for the global objective 
$\fone(\cdot)$, while the individual local objectives 
$\ftwo_i(\cdot,\cdot), i \in [N],$ may possibly be non-convex or non-smooth. 
Additionally, note that no differentiability is assumed for $\phi_i(\cdot)$.

We are now ready to present a convergence estimate for 
Algorithm~\ref{alg:ARGFree}, which is the main result of this paper. 
To this end, we first introduce the following instrumental notation:
\begin{align}\label{eq:theta_definition}
\theta_k&:=\operatorname{col}(\left\|x_k-x^*\right\|,
\left\|\sigma_k-\mathcal{J} \sigma_k\right\|,\left\|s_k-\mathcal{J} s_k\right\|,\\
&\hspace{3.5cm}\left\|\fp_k-\mathcal{J} \fp_k\right\|,
\left\|p_k-\mathcal{J} p_k\right\|), \notag
\end{align}
where $x^* \in \real^n$ is the unique (by~\ref{as:objective_a})
optimizer of~\eqref{eq:optimization_1}.

\begin{theorem}[\textbf{\textit{Convergence estimate for 
Algorithm~\ref{alg:ARGFree}}}]
\label{thm:algorithm_convergence}
Let Assumptions~\ref{as:graph}-\ref{as:objective} hold, suppose 
$\delta < \alpha \sqrt{n},$ 
$\Ltwo_0<\frac{1-\rho_A}{\|A-I\|}$,
and that the stepsize $\alpha$ satisfies
\begin{align} \label{eq:alpha_bound}
0 < \alpha < \min \left\{\frac{1}{2(n+4) \Lone_1}, \alpha_{1}^*,\alpha_{2}^* \right\},
\end{align}
where
\begin{align*}
\alpha^*_{1}&:=
\frac{1-\rho_A}{\Lphi\left(\sqrt{2(n+4)} \Lone_1+\frac{2 \sqrt{n}}{\delta}\right)},
\\
\alpha^*_{2}&:=
\frac{1-\rho_A-\Ltwo_0\|A-I\|}{\Ltwo_0\left(1+\Lphi\right)\left(\frac{2 \sqrt{n}}{\delta}+\sqrt{2(n+4)} \Lone_1\right)}.
\end{align*}
Then, there exists $\eta \in (0,1)$ and $\varepsilon \in \realpos$ such that,
for all $k \in \integnneg$, \eqref{eq:alg_vector_form} satisfies:
\begin{align}\label{eq_theta_bound} %
E[\theta_k]  \leq \eta^{k}
E[\theta_0]  + \frac{1-\eta^k}{1-\eta} \varepsilon.
\end{align}
Moreover\footnote{
We say that $\fone(x)=\mathcal{O}(g(x))$  as $x \rightarrow \infty$ if there 
exist constants $C>0$ and $x_0 \in \mathbb{R}$ such that:
$|\fone(x)| \leq C \cdot|g(x)|$ for all $x \geq x_0$.}, for large\footnote{That is, there exists $n_\circ,$ such that, for any $n \geq n_\circ$, the estimate holds.} $n,$
\begin{align}\label{eq:epsilon_order}
\varepsilon = 
\mc O( \delta \cdot E[\left\|u_{k+1}-u_k\right\|^2] \}) 
= \mathcal{O}\left(\delta  n \right).
\end{align}
% $$\varepsilon = \mathcal{O}\left(\delta^2 \cdot \max \left\{\alpha^2 n^3, \sup_{ t \in \integnneg, t \leq k-1 } E[\left\|u_{t+1}-u_t\right\|^2]\right\}\right).$$
\QEDB \end{theorem}

To maintain clarity in the exposition, the proof of this result is presented in 
Section~\ref{sec:analysis}.
Theorem~\ref{thm:algorithm_convergence} establishes that the 
iterates of Algorithm~\ref{alg:ARGFree} converge at rate~$\eta$ to a 
neighborhood of the  solution of~\eqref{eq:optimization_1}. 
Note that precise estimates for $\eta$ and $\epsilon$ are given shortly below 
(see Theorem~\ref{cor:exact_estimates}).
In analogy with~\eqref{eq:contraction_forward_difference}, convergence 
to an inexact point arises from employing a randomized gradient 
estimate instead of the exact gradient.
The algorithm’s accuracy $\varepsilon$ depends on the smoothing ratio 
$\delta$ and on the problem dimension~$n$, and can be controlled 
% (i.e., made arbitrarily small) 
by reducing $\delta.$
The algorithm's stepsize $\alpha$ is required to be sufficiently small, as 
given by the estimate~\eqref{eq:alpha_bound}. 
Note that the upper bound \( \alpha < \frac{1}{2(n+4)\Lone_1} \) recovers the 
maximum stepsize allowed for the centralized method (see 
Lemma~\ref{lem:bound_twopoint_alg}).
The additional bounds $\alpha < \alpha_1^*$ and $\alpha < \alpha_2^*$ are 
required to ensure contraction of the dynamic consensus tracking variables 
(lines 3--6 of Algorithm~\ref{alg:ARGFree}). 
Notice also that $\alpha_2^*$ is guaranteed to be a real positive number under 
the requirement $\Ltwo_0<\frac{1-\rho_A}{|A-I\|},$ which can be 
interpreted as an upper bound on the largest admitted variability on the 
functions $\ftwo_i$ (cf.~\ref{as:objective_c}) in relation to the 
averaging rate of the communication graph (measured by the parameters $\rho_A$ 
and $\norm{A-I}$).
Intuitively, the larger $\Ltwo_0,$ the faster the communication graph needs to 
be at averaging (cf. Assumption~\ref{as:graph}).
% As an illustrative example, suppose $A = (1-\lambda ) I + \lambda J$ so that
% $\rho_A = 1- \lambda$. As $\lambda \to 0,$ we have $\rho_A \to 1$ and $\norm{A-I}\to 0;$ as $\lambda \to 1,$ we have $\rho_A \to 0$ and $\norm{A-I}\to 2.$ 
We conclude by giving precise estimates for $\eta$ and $\varepsilon$ next.

\begin{theorem}[\textbf{\textit{Estimates for the convergence rate and 
accuracy of Algorithm~\ref{alg:ARGFree}}}]
\label{cor:exact_estimates}
Under the assumptions of Theorem~\ref{thm:algorithm_convergence}, the 
convergence rate in~\eqref{eq_theta_bound} satisfies:
\begin{align}\label{eq:bound_eta_star}
    \eta \leq \max \{ \eta_1^*, \eta_2^*, \eta_3^*\},
\end{align}
where
\begin{align*}
\eta_1^* &:= 
\sqrt{1 
-\alpha \mu\left(1-2 \alpha(n+4) \Lone_1\right)} + \frac{2 \alpha \sqrt{n}}{\delta},\\
\eta_2^* &:= \rho_A+\alpha \Lphi\left(\sqrt{2(n+4)} \Lone_1+\frac{2 \sqrt{n}}{\delta}\right),\\
\eta_3^* &:= 
\rho_A 
+ \alpha \Ltwo_0\left(1+\Lphi\right) \left(\sqrt{2(n+4)} \Lone_1 + \frac{2 \sqrt{n}}{\delta}\right)\\
& \quad\quad\quad +\Ltwo_0\|A-I\|.
\end{align*} 
Moreover, the approximation error $\varepsilon$ satisfies:
\begin{align}\label{eq:epsilon}
\varepsilon =  \delta \Bigg[& 
\alpha  \Lone_1 n
+ 2 \sbs{L}{comb}^2 E[\left\|u_{k+1}-u_k\right\|^2]^2 \\
& 
\quad\quad\quad\quad 
+ \frac{\alpha^2 \Lone_1^2}{2}(n+6)^3 \left(1+ \frac{3}{2} \sbs{L}{comb}^2\right) \Bigg]^{1/2},\notag
\end{align}
where $\sbs{L}{comb}^2 := \Lphi^2+\Ltwo_0^2\left(1+\Lphi\right)^2.$
\QEDB\end{theorem}

The proof of this result is presented in 
Section~\ref{sec:proof_cor_exact_estimates}. It is worth comparing the 
convergence rate estimate established by the Theorem~\ref{cor:exact_estimates}:
\[
\eta_1^* = \sqrt{1 - \alpha\mu \left(1 - 2 \alpha (n+4) \Lone_1 \right)} + \frac{2 n \alpha^2}{\delta^2},
\]
with the corresponding estimate for the centralized 
algorithm given in~\eqref{eq:contraction_forward_difference} (see Lemma~\ref{lem:bound_twopoint_alg}):
\[
\sqrt{1 - \beta_1^\alpha} = \sqrt{1 - \alpha \mu \left(1 - 2 \alpha (n+4)\Lone_1 \right)}.
\]
The degradation in rate affecting the distributed algorithm (characterized 
by the additional term $\frac{2 n \alpha^2}{\delta^2}$) can be traced 
back (as shown in the proof of 
Theorem~\ref{thm:algorithm_convergence}) to the presence of the 
tracking variables $\sigma_k, s_k, \fp_k, p_k.$ 
Intuitively, while the centralized algorithm has direct access to the 
quantities 
$
\sbs{\sigma}{f}(x_k), \quad \sbs{\sigma}{f}(x_k + \delta u_k), \quad \ftwo(x_k, \sigma_k), \quad \text{and} \quad \ftwo(x_k + \delta u_k, \sigma_k),
$
the distributed algorithm must rely on approximations of these quantities via 
the auxiliary variables $\sigma_k, s_k, \fp_k, p_k$, which in turn 
slows down the descent step.
The remaining estimates, $\eta_2^*$ and $\eta_3^*$, can be interpreted 
as bounds on the convergence rates of the tracking 
variables---specifically, $\eta_2^*$ corresponds to the convergence of 
$(\sigma_k, s_k)$, while $\eta_3^*$ pertains to $(\fp_k, p_k)$.

\begin{algorithm}[t]
\DontPrintSemicolon
\caption{\texttt{ARGFree-EM} (agent $i$)}
\label{alg:ARGFree_filtration}
\KwData{Parameters $\alpha, \delta \in \realpos,$ $B_i \in \real^{n_i \times n_i}, \Sigma_{u,i}^0, \Sigma_{v,i}\in \symmetric^{n_i}$}
\textbf{Initializations:} Set
$k=0,  u_{i}^{0} \sim \mc N(0, \Sigma_{u,i}^0)$, $x_{i}^{0} \in \real^{n_i},$ 
\begin{align*}
\sigma_{i}^{0} &= \phi_i(x_{i}^{0}), &
s_{i}^{0} &= \phi_i(x_{i}^{0} + \delta u_{i}^{0}), \\
\fp_{i}^{0} &= \ftwo_i(x_{i}^{0}, \sigma_{i}^{0}) &
p_{i}^{0} &= \ftwo_i(x_{i}^{0} + \delta u_{i}^{0}, s_{i}^{0})
\end{align*}

%
%\nl At each time $k\in \naturalnneg:$\\
%
\textbf{Optimization variable update:}\;
\nl
\begin{algomathdisplay}
x_i^{k+1} = x_i^k - \alpha \frac{p_i^k - \fp_i^k}{\delta} (\Sigma_{u,i}^k)^\inv u_i^k
\end{algomathdisplay}\;
\textbf{Exploration variables update:}\;
\nl Generate $v_{i}^{k+1} \sim \mc N(0, \Sigma_{v,i})$, 
$u_i^{k+1} = B_i u_i^k + v_i^{k+1}$\;
\nl Exploration covariance update:
\begin{algomathdisplay}
\Sigma_{u,i}^{k+1} = B_i \Sigma_{u,i}^{k} B_i^{\top}+\Sigma_{v,i}
\end{algomathdisplay}
\textbf{Tracking variables update:}\;
\begin{algomathdisplay}
\begin{aligned}
\nl \sigma_i^{k+1} &= \sum_{j\in \mc N_i} a_{ij} \sigma_j^k + \phi_i(x_i^{k+1})
- \phi_i(x_i^{k})\;\\
\nl s_i^{k+1} &= \sum_{j\in \mc N_i} a_{ij} s_j^k + \phi_i(x_i^{k+1} + \delta u_i^{k+1}) - \phi_i(x_i^{k} + \delta u_i^{k})\;\\
\nl \fp_i^{k+1} &= \sum_{j\in \mc N_i} a_{ij} \fp_j^k + \ftwo_i(x_i^{k+1},\sigma_i^{k+1}) - \ftwo_i(x_i^{k}, \sigma_i^k)\;\\
\nl p_i^{k+1} &= \sum_{j\in \mc N_i} a_{ij} p_j^k + \ftwo_i(x_i^{k+1} + \delta u_i^{k+1}, s_i^{k+1}) \\
& \hspace{3cm} - \ftwo_i(x_i^{k}+\delta u_i^{k}, s_i^k)\;\\
\end{aligned}
\end{algomathdisplay}\;
\nl Transmit $\sigma_i^{k+1}$, $s_i^{k+1}$, $\fp_i^{k+1}$ $p_i^{k+1}$ to neighbors\;
\nl Set $k \gets k+1$ and go to line 1\;
\KwResult{$x_i^k$, estimate for the optimizer of~\eqref{eq:optimization_1}}
\end{algorithm}

%%%%%%%%%%%%%%%%%%%%%%%%%%%%%%%%%%%%%%%%%%%%%%%%%%%%%%%%%%%%%%%%%%%%%%
%%%%%%%%%%%%%%%%%%%%%%%%%%%%%%%%%%%%%%%%%%%%%%%%%%%%%%%%%%%%%%%%%%%%%%%%%%%%%%%
\section{The \texttt{ARGFree-EM} algorithm: enhancing accuracy through exploration momentum}
\label{sec:single_realization}
While Theorem~\ref{thm:algorithm_convergence} guarantees convergence of 
Algorithm~\ref{alg:ARGFree} to a neighborhood of the optimizer 
of~\eqref{eq:optimization_1}, the size of this neighborhood can be large when 
$n$ is large. In this section, we propose a modification to \texttt{ARGFree} 
with improved accuracy.

%%%%%%%%%%%%%%%%%%%%%%%%%%%%%%%%%%%%%%%%%%%%%%%%%%%%%%%%%%%%%%%%%%%%%%
%%%%%%%%%%%%%%%%%%%%%%%%%%%%%%%%%%%%%%%%%%%%%%%%%%%%%%%%%%%%%%%%%%%%%%%%%%%%%%%
% \section{Improving the accuracy of \texttt{ARGFree} via exploration momentum}
\subsection{The \texttt{ARGFree-EM} algorithm}

Although Theorem~\ref{thm:algorithm_convergence} guarantees convergence of 
Algorithm~\ref{alg:ARGFree} to a neighborhood of the optimizer 
of~\eqref{eq:optimization_1}, the size of this neighborhood is dominated by 
\(\mathbb{E}[\|u_{k+1} - u_k\|^2]\), which can be large when the sequence 
$\{ u_k\}$ varies rapidly. More precisely, the estimate 
\(\mathbb{E}[\|u_{k+1} - u_k\|^2] = \mathcal{O}(n)\) (used 
to obtain~\eqref{eq:epsilon_order}) entails that the size of the convergence 
neighborhood (i.e., $\varepsilon$) scales with the problem dimension \(n\), thus potentially limiting performance in high-dimensional settings.

To overcome this limitation, a modified method for solving
Problem~\ref{prob:main_problem} is presented in 
Algorithm~\ref{alg:ARGFree_filtration}. The key innovation relative to 
Algorithm~\ref{alg:ARGFree} lies in the use of a filtered exploration signal 
(see line 2 of Algorithm~\ref{alg:ARGFree_filtration}):
\begin{align}
\label{eq:filtration}
u_i^{k+1} = B_i u_i^k + v_i^{k+1}, && i \in [N],
\end{align}
where $B_i \in \mathbb{R}^{n_i \times n_i}$ is a matrix encoding the momentum (or 
memory) of the filter, with spectral radius $\rho(B_i) < 1$, and 
$v_i^{k+1} \sim \mathcal{N}(0, \Sigma_{v,i})$ is a stochastic perturbation to 
the exploration signal,  with covariance 
$\Sigma_{v,i} \in \symmetric^{n_i}$. 
The filter is initialized at $u_i^0 \sim \mathcal{N}(0, \Sigma_{u,i}^0)$,
$\Sigma_{u,i}^0\in\symmetric^{n_i}$.
Because the filtration~\eqref{eq:filtration} introduces a correlation between 
sequential samples of the exploration signal $u_i^k$, the descent update (line 1 
of Algorithm~\ref{alg:ARGFree_filtration}) is modified (relative to 
Algorithm~\ref{alg:ARGFree}) with a normalization by the inverse correlation 
operator $(\Sigma_{u,i}^k)^\inv$, which is updated iteratively (see line 3 of 
Algorithm~\ref{alg:ARGFree_filtration}).

\subsection{Convergence guarantees for \texttt{ARGFree-EM}}

We next present a convergence estimate for 
Algorithm~\ref{alg:ARGFree_filtration}.
To state our result, we use the notation 
$B=\diag(B_1, \dots, B_N)\in \real^{n \times n}$ and
$\Sigma_v=\diag(\Sigma_{v,1}, \dots, \Sigma_{v,N}) \in \symmetric^n.$

\begin{theorem}[\textbf{\textit{Convergence estimate for 
Algorithm~\ref{alg:ARGFree_filtration}}}]
\label{thm:algorithm_convergence_filtration}
Let the assumptions of Theorem~\ref{thm:algorithm_convergence} hold.
Then, for large\footnote{Precisely, there exists $n_\circ,$ such that, for any $n \geq n_\circ$, the estimate holds.} $n,$ the iterates of Algorithm~\ref{alg:ARGFree_filtration} 
satisfy~\eqref{eq_theta_bound} with 
\begin{align}\label{eq:epsilon_order_filtration}
\varepsilon = 
 \mc O \left( \delta \cdot \max \{ \trace[(B-I)^{\top} \Sigma (B-I) + \Sigma_{v}]),  \sqrt{n}\}\right),
\end{align}
where $\Sigma \in \symmetric^{n}$ denotes the steady-state covariance of the 
process $\{u_k\}_{k \in \naturalnneg}$, given by the solution to the 
(implicit) equation $\Sigma = B \Sigma B^\top + \Sigma_v.$
\QEDB \end{theorem}

For ease of presentation, the proof of this result is postponed to 
Section~\ref{sec:pfof_algorithm_convergence_filtration}.
The bound~\eqref{eq:epsilon_order_filtration} shows that, incorporating 
momentum into the exploration signal (through the filtration 
model~\eqref{eq:filtration}) enables us to improve the accuracy of the 
algorithm. In particular, relative to Algorithm~\ref{alg:ARGFree}, we 
observe a reduction in the size of the convergence neighborhood by a worst-case 
factor of \(\sqrt{n}\): from 
\(\mathcal{O}(\delta n)\) in~\eqref{eq:epsilon_order} to 
\(\mathcal{O}(\delta \sqrt{n})\) in~\eqref{eq:epsilon_order_filtration}. 
Interestingly, \eqref{eq:epsilon_order_filtration} reveals the possibility of 
controlling the size of the convergence neighborhood (when the active term in 
the maximization of~\eqref{eq:epsilon_order_filtration} is 
$\trace[(B-I)^{\top} \Sigma (B-I) + \Sigma_{v}])$) by suitably selecting
the filtration parameters \(B_i\) and \(\Sigma_{v,i}\), as illustrated in Example~\ref{ex:choice_covariance}. 
The variables $B_i$ and $\Sigma_{v,i}$ thus serve as additional tuning 
parameters, providing extra flexibility beyond the smoothing ratio \(\delta\).

\begin{example}[\textbf{\textit{Choice of $B_i$ and $\Sigma_{v,i}$}}]
\label{ex:choice_covariance}
For instance, choosing \(B_i = \kappa I\) with \(\kappa \in (0,1)\) and \(\Sigma_{v,i} = \sigma^2 I\) yields a steady-state covariance
$\Sigma = B \Sigma B^\top + \Sigma_v = \kappa^2 \Sigma + \sigma^2 I,$
which solves to
\[
\Sigma = \frac{\sigma^2}{1 - \kappa^2} I.
\]
Substituting into the expression in~\eqref{eq:epsilon_order_filtration}, 
we obtain:
\[
\operatorname{Tr} \left[ (B - I)^\top \Sigma (B - I) + \Sigma_v \right] 
= (\kappa - 1)^2  \frac{n \sigma^2}{1 - \kappa^2} + n \sigma^2.
\]
This quantity can be reduced by selecting sufficiently-small parameters 
\(\sigma\) and \(\kappa\). 
% \(\sigma, \kappa\).~
% The variables \(\sigma, \kappa\) thus serve as additional tuning 
% parameters, providing extra flexibility beyond the smoothing ratio \(\delta\).
\QEDB\end{example}

%%%%%%%%%%%%%%%%%%%%%%%%%%%%%%%%%%%%%%%%%%%%%%%%%%%%%%%%%%%%%%%%%%%%%%
%%%%%%%%%%%%%%%%%%%%%%%%%%%%%%%%%%%%%%%%%%%%%%%%%%%%%%%%%%%%%%%%%%%%%%%%%%%%%%%
\section{Convergence analysis of the algorithms}
\label{sec:analysis}
This section is devoted to establishing the bounds presented 
theorems~\ref{thm:algorithm_convergence}, \ref{cor:exact_estimates}, 
and~\ref{thm:algorithm_convergence_filtration}. Our approach is based on 
showing that, for all $k \in \integnneg$, the quantity $\theta_k$ 
satisfies the following bound:
\begin{align}\label{eq:recursion_theta}
E[\theta_{k+1}] \leq M(\alpha) E[\theta_k] + b,
\end{align}
where $M(\alpha)$ is a Schur stable matrix with spectral radius $\eta$ and
$b$ is such that $\norm{b}\leq \varepsilon.$ 
We begin with the proof of Theorem~\ref{thm:algorithm_convergence}, which is 
organized into seven subsections 
(Section~\ref{sec:bound_xkk}--\ref{sec:pfof_thm_algorithm_convergence}). 
We conclude by presenting the  proof of Theorem~\ref{cor:exact_estimates} in 
Section~\ref{sec:proof_cor_exact_estimates}
and that of Theorem~\ref{thm:algorithm_convergence_filtration} in 
Section~\ref{sec:pfof_algorithm_convergence_filtration}.

\subsection{Bound for $E[\left\|x_{k+1}-x^*\right\|]$}
\label{sec:bound_xkk}
We have the following estimate:
\begin{align*}
E_{u_k}&[\left\|x_{k+1}-x^*\right\|] 
= E_{u_k}[\norm{x_k -\frac{\alpha}{\delta} (p_k-\fp_k)u_k - x^*}]\\
&\leq \underbrace{E_{u_k}[ \norm{x_k-\alpha g_\delta(x_k) -x^*}]}_{:=\mathcircled{a}} \\
&\quad + \underbrace{\alpha E[\norm{g_\delta(x_k) - \frac{p_k-\fp_k}{\delta}u_k}]}_{:=\mathcircled{b}}.
\end{align*}
By application of Lemma~\ref{lem:bound_twopoint_alg}, we have 
$\mathcircled{a}\leq\sqrt{1-\beta_1^\alpha} \norm{x_k - x^*} + \beta_2^\alpha;$ 
The term $\mathcircled{b}$ satisfies:
\begin{align}\label{eq:bound_b}
\mathcircled{b} &\leq
\underbrace{\alpha E_{u_k}[\norm{g_\delta(x_k) - \one \otimes \frac{\bar p_k-\bar \fp_k}{\delta}u_k}]}_{:=\mathcircled{c}}\notag\\
&\quad + \underbrace{\frac{\alpha}{\delta} E_{u_k}[\norm{(p_k - \one \otimes \bar p_k)u_k}]}_{:=\mathcircled{d}}\notag\\
&\quad + \underbrace{\frac{\alpha}{\delta} E_{u_k}[\norm{(\fp_k - \one \otimes \bar \fp_k)u_k}]}_{:=\mathcircled{e}},
\end{align}
where the first inequality follows by adding and subtracting $\one \otimes \frac{\bar p_k-\bar \fp_k}{\delta}u_k$ inside the norm.
Next, observe that $ \mathcircled{c} \reljust{=}{\eqref{eq:averages_identities},\eqref{eq:grad_oracle}} 0;$ the remaining terms satisfy:
\begin{align}\label{eq:bound_d_e}
\mathcircled{d} &\reljust{\leq}{\eqref{eq:kronone_calJ},\eqref{eq:expectation_norm_bound}} 
\frac{\alpha \sqrt{n}}{\delta} \norm{p_k - \mc J p_k}, \notag\\
\mathcircled{e} &\reljust{\leq}{\eqref{eq:kronone_calJ},\eqref{eq:expectation_norm_bound}} 
\frac{\alpha \sqrt{n}}{\delta} \norm{\fp_k - \mc J \fp_k}.
\end{align}
Summarizing, we have derived the estimate:
\begin{align}\label{eq:xkp1mxstar}
E_{u_k}[\left\|x_{k+1}-x^*\right\|] &\leq 
\sqrt{1-\beta_1^\alpha} \norm{x_k - x^*} + \beta_2^\alpha\\
& \quad  + \gamma_1^{\alpha} (
\norm{\fp_k - \mc J \fp_k} +
\norm{p_k - \mc J p_k}), \notag
\end{align}
where we defined $\gamma_1^{\alpha} =: \alpha \sqrt{n} /\delta$.

\subsection{An auxiliary bound for $E[\left\|x_{k+1}-x_k\right\|]$}
We have the following estimate:
\begin{align}\label{eq:xkp1mxk}
E_{u_k}&[\left\|x_{k+1}-x_k\right\|] = \alpha E[\norm{\frac{p_k -\fp_k}{\delta}u_k}]\notag\\
& \leq \underbrace{\alpha E_{u_k}[\norm{g_\delta(x_k)-\frac{p_k-\fp_k}{\delta}u_k}]}_{=\mathcircled{b}} + \alpha E_{u_k}[\norm{g_\delta(x_k)}]\notag\\
&\hspace{-.3cm}\reljust{\leq}{\eqref{eq:bound_b},\eqref{eq:bound_d_e}}
\gamma_1^{\alpha} (\norm{p_k - \mc J p_k} + \norm{\fp_k - \mc J \fp_k})\notag\\
&\quad\quad+\alpha E_{u_k}[\norm{g_\delta(x_k)}]\notag\\
& \reljust{\leq}{\eqref{eq:Nesterov35}}
\gamma_1^{\alpha} (\norm{p_k - \mc J p_k} + \norm{\fp_k - \mc J \fp_k})\notag\\
 &\quad\quad\quad + \alpha \left(\frac{\delta^2}{2} \Lone_1^2(n+6)^3 + 2(n+4)\|\nabla \fone(x )\|^2 \right)^{1/2}\notag\\
&\leq \gamma_1^{\alpha} (\norm{p_k - \mc J p_k} + \norm{\fp_k - \mc J \fp_k}) + \gamma_2^{\alpha} \notag\\
&\quad + \gamma_3^{\alpha} \norm{x_k - x^*},
\end{align}
where for the third inequality we combined~\eqref{eq:Nesterov35} with Jensen's 
inequality (which, for a scalar $a\geq 0$, gives 
$E [\sqrt{a}]\leq \sqrt{E[a]}$) and the 
last inequality we used $\norm{\nabla \fone(x)} = \norm{\nabla \fone(x) -\nabla \fone(x^*)}\leq \Lone_1 \norm{x-x^*}$ and
we defined 
$\gamma_2^{\alpha}:=\frac{1}{2} \alpha \delta \Lone_1 (n+6)^{3/2}$ and 
$\gamma_3^{\alpha}:=\sqrt{2 (n+4)}\alpha \Lone_1.$
% Here, to derive the third line, we recognized the term 
% $\mathcircled{b}$ bounded in~\eqref{eq:bound_b}.

\subsection{Bound for $E[\left\|\sigma_{k+1}-\mc J \sigma_{k+1}\right\|]$}
We have the estimate:
\begin{align}\label{eq_sigmakp1msigmakp1}
E_{u_k}&[\norm{\sigma_{k+1}-\mc J \sigma_{k+1}}] \\
&\reljust{\leq}{\eqref{eq:alg_vector_form}}
\norm{\mc A \sigma_k - \mc J \mc A\sigma_k}\notag\\
&\quad\quad 
+ E_{u_k}[\norm{(I - \mc J)(\sbs{\phi}{v}\left(x_{k+1}\right)- \sbs{\phi}{v} \left(x_k\right))}]\notag\\
&
\reljust{\leq}{\eqref{eq:inequalities_AJ}}
\rho_A \norm{\sigma_k - \mc J\sigma_k}\notag\\
&\quad +  E_{u_k}[\norm{\sbs{\phi}{v}\left(x_{k+1}\right)- \sbs{\phi}{v} \left(x_k\right)}]
\notag\\
&\reljust{\leq}{\ref{as:objective_b}}
\rho_A \norm{\sigma_k - \mc J \sigma_k} + \Lphi E_{u_k}[\norm{x_{k+1}-x_k}]\notag\\
&\reljust{\leq}{\eqref{eq:xkp1mxk}}
\rho_A \norm{\sigma_k - \mc J \sigma_k} + \Lphi \gamma_1^{\alpha}
\norm{p_k - \mc J p_k} \notag\\
&\quad +  \Lphi \gamma_1^{\alpha} \norm{\fp_k - \mc J \fp_k} + \Lphi\gamma_2^{\alpha}  + \Lphi \gamma_3^{\alpha} \norm{x_k - x^*}, \notag
\end{align}
where for the second inequality we used $\norm{I-\mc J}=1.$

\subsection{Bound for $E[\left\|s_{k+1}-\mc J s_{k+1}\right\|]$}
Notice that $s_{k+1}$ has two sources of stochasticity at time $k$: $u_k$ and 
$u_{k+1}.$ We thus have:
\begin{align}\label{eq_skp1mskp1}
&\hspace{-.4cm} E_{u_k,u_{k+1}}[\left\|s_{k+1}-\mc J s_{k+1}\right\|]\\
&\reljust{\leq}{\eqref{eq:alg_vector_form}}
\norm{\mc A s_k - \mc J \mc As_k}\notag\\
&\quad\quad 
+ E_{u_k,u_{k+1}}[\norm{(I - \mc J)(\sbs{\phi}{v}\left(x_{k+1}+\delta u_{k+1}\right) \notag\\
&\hspace{5cm} - \sbs{\phi}{v} \left(x_k+\delta u_k\right))}]\notag\\
&\reljust{\leq}{\eqref{eq:inequalities_AJ}}
\rho_A \norm{ s_k - \mc J s_k}\notag\\
&\quad\quad +  E_{u_k,u_{k+1}}[\norm{\sbs{\phi}{v}\left(x_{k+1}+\delta u_{k+1}\right)- \sbs{\phi}{v} \left(x_k + \delta u_k\right)}]
\notag\\
&\reljust{\leq}{\ref{as:objective_b}}
\rho_A \norm{ s_k - \mc J s_k} + \Lphi E_{u_k}[\norm{x_{k+1}-x_k}]\notag\\
&\quad\quad  + \Lphi\delta E_{u_k,u_{k+1}}[\norm{u_{k+1}-u_k}]\notag\\
&\reljust{\leq}{\eqref{eq:xkp1mxk}}
\rho_A \norm{s_k - \mc J s_k} + \Lphi \gamma_1^{\alpha}
\norm{p_k - \mc J p_k} \notag\\
&\quad\quad +  \Lphi \gamma_1^{\alpha} \norm{\fp_k - \mc J \fp_k} + \Lphi\gamma_2^{\alpha}  + \Lphi \gamma_3^{\alpha} \norm{x_k - x^*}\notag\\
&\quad\quad  + \Lphi\delta E_{u_k,u_{k+1}}[\norm{u_{k+1}-u_k}], \notag
\end{align}
where for the second inequality we used $\norm{I-\mc J}=1.$

\subsection{Bound for $E[\left\|\fp_{k+1}-\mc J \fp_{k+1}\right\|]$}
Before bounding the desired term, notice that:
\begin{align}\label{eq_sigmakp1msigmak}
E_{u_k}&[\norm{\sigma_{k+1}-\sigma_{k}}] \\
&\reljust{\leq}{\eqref{eq:alg_vector_form}}
\norm{\mc A \sigma_k - \sigma_k} 
+ E_{u_k}[\norm{\sbs{\phi}{v}\left(x_{k+1}\right)- \sbs{\phi}{v} \left(x_k\right))}]\notag\\
&\reljust{\leq}{\ref{as:objective_b}}
\norm{(\mc A - I)(\sigma_k - \mc J \sigma_k)} 
+ \Lphi E_{u_k}[\norm{x_{k+1}-x_k}]\notag\\
&\quad\reljust{\leq}{\eqref{eq:xkp1mxk}}
\norm{A - I} \norm{\sigma_k - \mc J \sigma_k} + \Lphi \gamma_3^{\alpha} \norm{x_k - x^*} \notag\\
&\quad\quad
+\Lphi \gamma_1^{\alpha} (\norm{p_k - \mc J p_k} + \norm{\fp_k - \mc J \fp_k}) + \Lphi \gamma_2^{\alpha}, \notag
\end{align}
where, for  the second inequality, we used the identity 
$\mc A \sigma_k -\sigma_k 
\reljust{=}{\eqref{eq:inequalities_AJ}} (\mc A - I) (\sigma_k - \mc J \sigma_k).$
We then have:
\begin{align}\label{eq_varphikp1mvarphikp1}
E_{u_k}&[\left\|\fp_{k+1}-\mc J \fp_{k+1}\right\|] \\
&\reljust{\leq}{\eqref{eq:alg_vector_form}}
\norm{\mc A \fp_k - \mc J \mc A\fp_k}\notag\\
&\quad\quad 
+ E_{u_k}[\norm{(I - \mc J)(\sbs{\ftwo}{v}\left(x_{k+1},\sigma_{k+1}\right)- \sbs{\ftwo}{v} \left(x_k,\sigma_k\right))}]\notag\\
&\reljust{\leq}{\eqref{eq:inequalities_AJ},\ref{as:objective_c}}
\rho_A \norm{\fp_k - \mc J\fp_k} \notag \\
&\quad\quad + \Ltwo_0 (E_{u_k}[\left\|x_{k+1}-x_k\right\|] + 
E_{u_k}[\left\|\sigma_{k+1}-\sigma_k\right\|])\notag\\
&
\reljust{\leq}{\eqref{eq:xkp1mxk},\eqref{eq_sigmakp1msigmak}}
 \gamma_4^\alpha \norm{x_k - x^*} 
+ \gamma_5 \norm{\sigma_k - \mc J \sigma_k}\notag\\
&\quad\quad
+ (\rho_A +\gamma_6^\alpha) \norm{\fp_k - \mc J \fp_k}
+ \gamma_6^\alpha  \norm{p_k - \mc J p_k}\notag
+ \gamma_7^\alpha,
\end{align}
where in the last step we defined
$\gamma_4^\alpha := \Ltwo_0 \gamma_3^\alpha (1+\Lphi),$
$\gamma_5 := \Ltwo_0 \norm{A-I},$
$\gamma_6^\alpha :=  \Ltwo_0 \gamma_1^\alpha(1+\Lphi),$ and
$\gamma_7^\alpha := \Ltwo_0 \gamma_2^\alpha(1+\Lphi).$

\subsection{Bound for $E[\left\|p_{k+1}-\mc J p_{k+1}\right\|]$}
We first need the following intermediate bound:
\begin{align}\label{eq_skp1msk}
&\hspace{-.4cm}
E_{u_k,u_{k+1}}[\norm{s_{k+1}-s_{k}}] \\
& \reljust{\leq}{\eqref{eq:alg_vector_form}}
\norm{\mc A s_k - s_k} \notag\\
&\quad\quad+ E_{u_k}[\norm{\sbs{\phi}{v}\left(x_{k+1}+\delta u_{k+1}\right)- \sbs{\phi}{v} \left(x_k + \delta u_k\right))}]\notag\\
& \reljust{\leq}{\ref{as:objective_b}}
\norm{(\mc A - I)(s_k - \mc J s_k)} \notag\\
&\quad\quad  + \Lphi E_{u_k}[\norm{x_{k+1}-x_k}]
+ \Lphi \delta E_{u_k,u_{k+1}}[\norm{u_{k+1}-u_k}] \notag\\
&\quad\reljust{\leq}{\eqref{eq:xkp1mxk}}
\norm{A - I} \norm{s_k - \mc J s_k} + \Lphi \gamma_3^{\alpha} \norm{x_k - x^*} \notag\\
&\quad\quad
+\Lphi \gamma_1^{\alpha} (\norm{p_k - \mc J p_k} + \norm{\fp_k - \mc J \fp_k}) + \Lphi \gamma_2^{\alpha}, \notag\\
&\quad\quad  
+ \Lphi \delta E_{u_k,u_{k+1}}[\norm{u_{k+1}-u_k}]. \notag
\end{align}
Then, 
\begin{align}\label{eq_pkp1mpkp1}
& \hspace{-.4cm} 
E_{u_k,u_{k+1}}[\left\|p_{k+1}-\mc J p_{k+1}\right\|]\\
& \reljust{\leq}{\eqref{eq:alg_vector_form}}
\norm{\mc A p_k - \mc J \mc A p_k}\notag\\
&\quad\quad 
+ E_{u_k,u_{k+1}}[\norm{(I - \mc J)(\sbs{\ftwo}{v}\left(x_{k+1}+\delta u_{k+1},s_{k+1}\right) \notag\\
&\quad\quad\quad\quad\quad\quad\quad\quad\quad\quad\quad\quad\quad\quad - \sbs{\ftwo}{v} \left(x_k+\delta u_{k},s_k\right))}]\notag\\
&\reljust{\leq}{\eqref{eq:inequalities_AJ},\ref{as:objective_c}}
\rho_A \norm{p_k - \mc J p_k} 
+ 2 \Ltwo_0 E_{u_k}[\left\|x_{k+1}-x_k\right\|] \notag\\
& \quad\quad + \Ltwo_0 E_{u_k,u_{k+1}}[\left\|s_{k+1}-s_k\right\|]\notag\\
& \quad\quad + \Ltwo_0  \delta E_{u_k,u_{k+1}}[\left\|u_{k+1}-u_k\right\|] \notag\\
&
\reljust{\leq}{\eqref{eq:xkp1mxk},\eqref{eq_sigmakp1msigmak}}
 \gamma_4^\alpha \norm{x_k - x^*}
+ \gamma_5 \norm{s_k - \mc J s_k} \notag\\
&\quad\quad + \gamma_6^\alpha \norm{\fp_k - \mc J \fp_k}
+ (\rho_A+ \gamma_6^\alpha) \norm{p_k - \mc J p_k} \notag\\
&\quad\quad +\gamma_7^\alpha
+ \gamma_8 E_{u_k,u_{k+1}}[\norm{u_{k+1} - u_k}], \notag
\end{align}
where we defined $\gamma_8:= \Ltwo_0 \delta (1+\Lphi).$

\subsection{Proof of Theorem~\ref{thm:algorithm_convergence}}
\label{sec:pfof_thm_algorithm_convergence}

By combining the estimates \eqref{eq:xkp1mxstar}, \eqref{eq_sigmakp1msigmakp1},
\eqref{eq_skp1mskp1}, \eqref{eq_varphikp1mvarphikp1}, and \eqref{eq_pkp1mpkp1}, 
we conclude that~\eqref{eq:recursion_theta} holds with 
\begin{align}\label{eq:M_alpha_explicit}
M(\alpha) =\left[\begin{array}{ccccc}\sqrt{1-\beta_1^\alpha} & 0 & 0 & \gamma_1^\alpha & \gamma_1^\alpha \\ 
\Lphi \gamma_3^\alpha & \rho_A & 0 & \Lphi \gamma_1^\alpha & \Lphi \gamma_1^\alpha \\ \Lphi \gamma_3^\alpha & 0 & \rho_A & \Lphi \gamma_1^\alpha & \Lphi \gamma_1^\alpha \\ \gamma_4^\alpha & \gamma_5 & 0 & \rho_A+\gamma_6^\alpha & \gamma_6^\alpha \\ \gamma_4^\alpha & 0 & \gamma_5 & \gamma_6^\alpha & \rho_A + \gamma_6^\alpha\end{array}\right]
\end{align}
and
\begin{align*}
b = \begin{bmatrix}
\beta_2^\alpha \\
\Lphi \gamma_2^\alpha \\
\Lphi \gamma_2^\alpha   + \Lphi \delta E_{u_k,u_{k+1}}[\left\|u_{k+1}-u_k\right\|] \\
\gamma_7^\alpha \\
\gamma_7^\alpha + \gamma_8 E_{u_k,u_{k+1}}[\left\|u_{k+1}-u_k\right\|]
\end{bmatrix}.
\end{align*}
Notice that $b$ satisfies 
$\norm{b} \leq \mc O( \delta \cdot E[\left\|u_{k+1}-u_k\right\|^2] \}).$
By substituting the expressions of the constants involved, we have that, 
for\footnote{We say $f(x)=o(g(x))$ as $x \rightarrow a$ if 
$\lim _{x \rightarrow a} \frac{f(x)}{g(x)}=0$.} $\alpha \to 0^+,$
\begin{align*}
M(\alpha) &= \begin{bmatrix} 
\sqrt{1 -\alpha \mu}& 0 & 0 & 0 &0\\
0 & \rho_A & 0 & 0 &0\\
0 &0 & \rho_A & 0 & 0 \\
0 & \gamma_5 & 0 & \rho_A & 0\\
0 & 0 & \gamma_5 &0 & \rho_A \\
\end{bmatrix}
+ o(\alpha),
\end{align*}
which proves that the spectral radius of $M(\alpha)$ is smaller than $1$ 
for sufficiently small $\alpha>0$.
To determine an estimate on the largest value of $\alpha$ that guarantees  $
\rho(M(\alpha))<1,$ we leverage the Gershgorin Circle Theorem. An application 
of the theorem to row $1$ of $M(\alpha)$ gives the condition:
\begin{subequations}\label{eq:gershgorin}
\begin{align}\label{eq:gershgorin_a}
\sqrt{1 
-\alpha \mu\left(1-2 \alpha(n+4) \Lone_1\right)} <1 - \frac{2 \alpha \sqrt{n}}{\delta}.
\end{align}
Squaring both sides of~\eqref{eq:gershgorin_a} and simplifying yields:
$$
\mu(1 - 2\alpha(n+4)\Lone_1) > \frac{4\sqrt{n}}{\delta} - \frac{4\alpha n}{\delta^2}.
$$
The condition $\alpha < \frac{1}{2(n+4)\Lone_1}$ ensures that the left-hand 
side is positive, while the assumption $\delta \leq \alpha \sqrt{n}$ 
guarantees that the right-hand side is non-negative. Hence, under the 
assumptions of Theorem~\ref{thm:algorithm_convergence}, the above inequality 
holds, and consequently~\eqref{eq:gershgorin_a} is also satisfied.
Applying Gershgorin's theorem to rows 2--3 of $M(\alpha)$ gives:
\begin{align}\label{eq:gershgorin_b}
\rho_A+\alpha \Lphi\left(\sqrt{2(n+4)} \Lone_1+\frac{2 \sqrt{n}}{\delta}\right)<1
% \rho_A^2+\alpha^2 \Lphi^2\left(2(n+4) L_1^2+\frac{2 n}{\delta^2}\right)<1,
\end{align}
which gives the estimate $\alpha_1^*$; applying Gershgorin's theorem to rows 
4-5 of $M(\alpha)$ gives:
\begin{align}\label{eq:gershgorin_c}
& \rho_A +\Ltwo_0\|A-I\|\\
& \quad\quad\quad +\Ltwo_0\left(1+\Lphi\right) \alpha \left(\frac{2 \sqrt{n}}{\delta} +\sqrt{2(n+4)} \Lone_1\right)<1 \notag
\end{align}
\end{subequations}
which gives the estimate $\alpha_2^*.$

We have thus shown that the estimate~\eqref{eq_theta_bound} holds, 
where \(M(\alpha)\) is a Schur-stable matrix with spectral radius 
\(\eta = \rho(M(\alpha))\), and \(b\) is a vector satisfying 
\(\|b\| = \mathcal{O}(\delta \, \mathbb{E}[\|u_{k+1}-u_k\|^2])\). 
Hence, the bound~\eqref{eq_theta_bound} is established with 
\(\varepsilon = \mathcal{O}(\delta \, \mathbb{E}[\|u_{k+1}-u_k\|^2])\), 
consistent with the first identity in~\eqref{eq:epsilon_order}.

Finally, we are left to prove the last identity in~\eqref{eq:epsilon_order}; 
that is, $\norm{b} = \mathcal{O}\left(\delta  {n} \right).$ To this end, 
we apply the inequality $(a + b)^2 \leq 2a^2 + 2b^2$ and bound the norm of the 
vector $b$ as:
\begin{align}\label{eq:bound_norm_b}
\|b\| \leq \delta \Bigg[& 
\alpha  \Lone_1 n+ \frac{\alpha^2 \Lone_1^2}{2}(n+6)^3 \\
& + 2 \left(\Lphi^2+\Ltwo_0^2\left(1+\Lphi\right)^2\right) E[\left\|u_{k+1}-u_k\right\|^2] \notag\\
& + \frac{3}{4} \alpha^2 \Lone_1^2(n+6)^3\left(\Lphi^2+\Ltwo_0^2\left(1+\Lphi\right)^2\right) \Bigg]^{1/2} \notag\\
\leq \delta \Bigg[& 
 \alpha \Lone_1 n 
 +8 n^2\left(\Lphi^2+\Ltwo_0^2\left(1+\Lphi\right)^2\right) \notag\\
 & +\alpha^2 \Lone_1^2(n+6)^3\left(\frac{1}{2}+\frac{3}{4} \Lphi^2+\frac{3}{4} \Ltwo_0^2\left(1+\Lphi\right)^2\right) 
\Bigg]^{1/2} \notag\\
&= \delta \mc O(\max \{ \sqrt{\alpha n}, n, \alpha n^{3/2}\}) \notag.
\end{align}
where the last inequality follows from
$E[\left\|u_{k+1}-u_k\right\|^2]\leq 2n$ (which follows 
from~\eqref{eq:expectation_norm_bound}, since  $u_{k+1}$ and $u_{k}$ are 
independent). By noting that $\alpha =\mc O(\frac{1}{n}),$ it follows that 
the second term dominates the maximization for large $n$ and the estimate 
$\norm{b} \leq \varepsilon = \mathcal{O}\left(\delta  n \right)$ follows.

\subsection{Proof of Theorem~\ref{cor:exact_estimates}}
\label{sec:proof_cor_exact_estimates}
The estimate~\eqref{eq:bound_eta_star} follows by noting that, 
from~\eqref{eq:recursion_theta}, $\eta$ can be selected to be an upper bound 
for the spectral radius of $M(\alpha)$ and, by 
using~\eqref{eq:M_alpha_explicit} and~\eqref{eq:gershgorin}, 
\eqref{eq:bound_eta_star} follows. 
Finally, the expression~\eqref{eq:epsilon} follows
from~\eqref{eq:recursion_theta} and  the estimate \eqref{eq:bound_b}.

\subsection{Proof of Theorem~\ref{thm:algorithm_convergence_filtration}}
\label{sec:pfof_algorithm_convergence_filtration}

We begin by noting that, under the filtration model~\eqref{eq:filtration}, the covariance of the process \((u_{k+1} - u_k)\) is:
\begin{align*}
\mathbb{E}\left[\|u_{k+1} - u_k\|^2\right] 
&= \mathbb{E}\left[\|N u_k + v_{k+1} - u_k\|^2\right] \notag\\
&= \mathbb{E}\left[\|(N - I) u_k + v_{k+1}\|^2\right] \notag\\
&= \mathbb{E}\left[u_k^\top (N - I)^\top (N - I) u_k\right] \notag \\
 &\quad\quad\quad + \mathbb{E}\left[v_{k+1}^\top v_{k+1}\right],
\end{align*}
where we used the independence and zero-mean property of \(v_{k+1}\) to eliminate the cross term. This yields that, in the stationary regime, the covariance of the process \((u_{k+1} - u_k)\) is given by
\begin{align}\label{covariance_derivation_u}
\operatorname{Tr}\left[(N - I)^\top \Sigma (N - I) + \Sigma_v\right],
\end{align}
where \(\Sigma\) is the stationary covariance matrix of \(u_k\), and \(\Sigma_v\) is the covariance of the noise term \(v_k\).
Next, by iterating through the steps of the proof of 
Theorem~\ref{thm:algorithm_convergence}, we obtain,  from~\eqref{eq:bound_norm_b}:
\begin{align}\label{eq:bound_norm_b_filtration}
\|b\| &
\reljust{\leq}{\eqref{eq:bound_norm_b}}
\delta \Bigg[ 
\alpha  \Lone_1 n+ \frac{\alpha^2 \Lone_1^2}{2}(n+6)^3 \\
& + 2 \left(\Lphi^2+\Ltwo_0^2\left(1+\Lphi\right)^2\right) E[\left\|u_{k+1}-u_k\right\|^2]^2 \notag\\
& + \frac{3}{4} \alpha^2 \Lone_1^2(n+6)^3\left(\Lphi^2+\Ltwo_0^2\left(1+\Lphi\right)^2\right) \Bigg]^{1/2} \notag\\
&\reljust{=}{\eqref{covariance_derivation_u}}
\mc O(\delta  \cdot \max \{   \sqrt{\alpha n},  \alpha n^{3/2},\notag\\
&\quad\quad\quad\quad\quad   \trace[(N-I)^{\top}(N-I) \Sigma + \Sigma_{v}]\}) \notag\\
&= \mc O(\delta \cdot \max \{ \sqrt{n}, \notag\\
&\quad\quad\quad\quad\quad \trace[(N-I)^{\top}(N-I) \Sigma + \Sigma_{v}]\}), \notag
\end{align}
where the last identity follows from the observation that, since 
\(\alpha = \mathcal{O}\left(\frac{1}{n}\right)\), the second term dominates the 
first one. This concludes the proof.

\begin{figure*}[t]
\centering \subfigure[]{\includegraphics[width=.63\columnwidth]{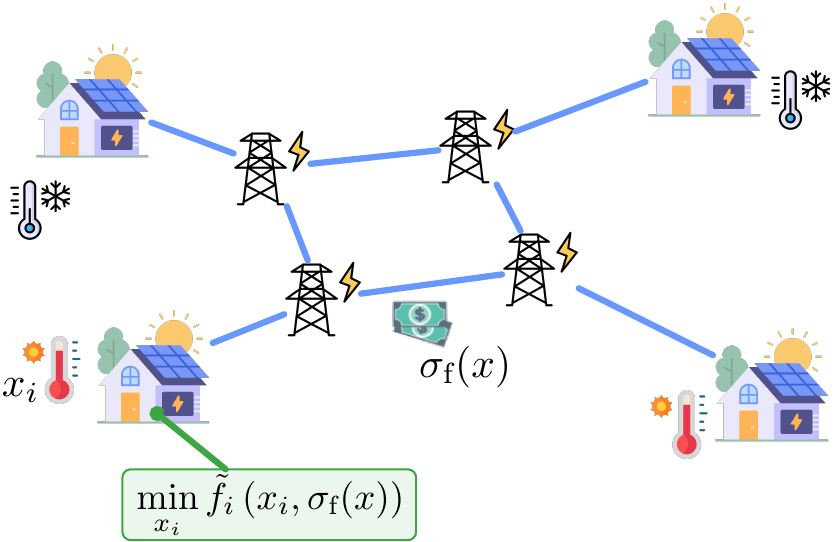}}\hfill
\centering \subfigure[]{\includegraphics[width=.63\columnwidth]{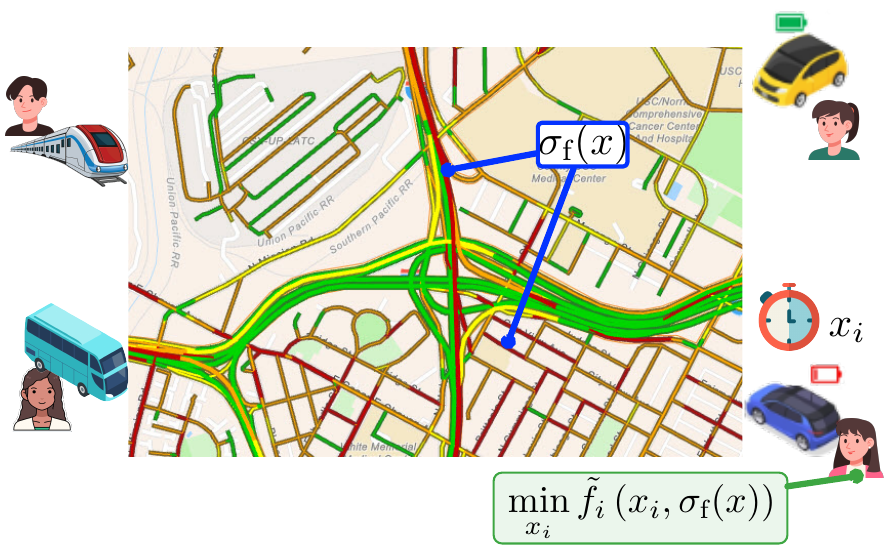}} \hfill
\centering \subfigure[]{\includegraphics[width=.63\columnwidth]{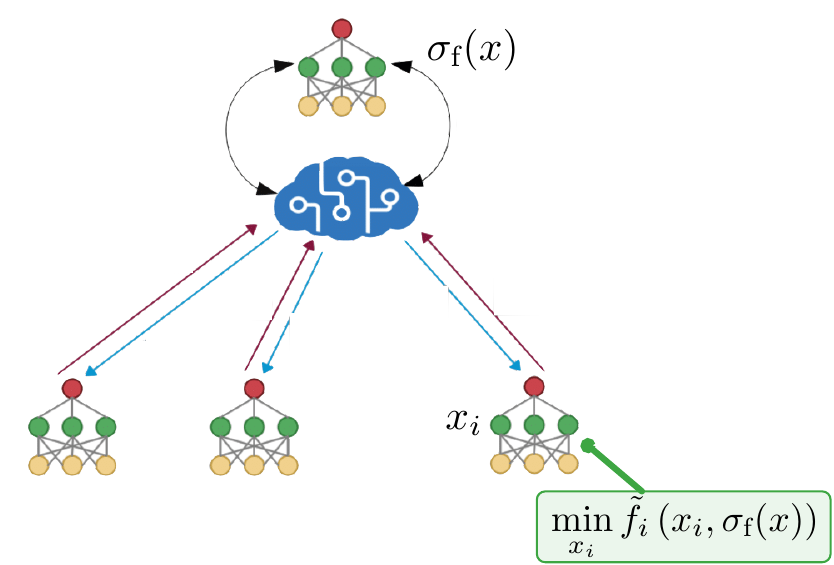}}
\caption{
Illustrative applications of the aggregative cooperative optimization 
framework~\eqref{eq:sigma_x}--\eqref{eq:optimization_1}, 
subject to requirements~\ref{requirement_1}--\ref{requirement_2}. 
(a) Application to energy sharing markets, where a prosumer seeks to 
decide how much energy to consume/inject to the grid while accounting 
for variable electricity prices. 
(b) application to transportation systems, where a user seeks to 
select a mode of transportation and experiences congestion based on 
collective transportation choices. 
(c) application to federated learning, where devices seek to optimize 
the prediction performance of a centralized machine learning model, 
while maintaining local data privacy. See
Section~\ref{sec:illustrative_applications} for a detailed discussion.
}
\vspace{-.2cm}
\label{fig:applications}
\end{figure*}

%%%%%%%%%%%%%%%%%%%%%%%%%%%%%%%%%%%%%%%%%%%%%%%%%%%%%%%%%%%%%%%%%%%%%%
%%%%%%%%%%%%%%%%%%%%%%%%%%%%%%%%%%%%%%%%%%%%%%%%%%%%%%%%%%%%%%%%%%%%%%%%%%%%%%%
\section{Illustrative applications}
\label{sec:illustrative_applications}
We next present four representative applications where 
problem~\eqref{eq:sigma_x}--\eqref{eq:optimization_1} and 
requirements~\ref{requirement_1}--\ref{requirement_2} arise naturally.

\textit{(Multi-agent robotics)}
Consider a multi-agent robotic coordination task (see 
Fig.~\ref{fig:applications_robotics}), where each agent \( i \) represents a 
robot that controls its position \( x_i \in \mathbb{R}^d \). 
The group aims to achieve a desired formation, for which a key quantity is the 
barycenter of the team, defined as
$\sbs{\sigma}{f}(x) = \frac{1}{N} \sum_{i=1}^N x_i.$
The local costs in~\eqref{eq:optimization_1} can then be used to reflect local 
positioning objective and a coupling to the group’s collective motion:
\begin{align}\label{eq:example-robotics}
\ftwo_i(x_i, \sbs{\sigma}{f}(x)) &= \text{local positioning cost}(x_i) \notag  \\
& \quad\quad + \text{formation cost}(\|x_i - \sbs{\sigma}{f}(x)\|).
\end{align}
In practice, the absolute position $x_i$ is often difficult to estimate 
reliably. Agents typically have access only to a function of $x_i$ (e.g., via 
cameras or perception systems), which itself is hard to determine accurately 
due  to the need for accurate calibration. In such cases, 
\eqref{eq:example-robotics} 
becomes:
\begin{align*}
\ftwo_i(x_i, \sbs{\sigma}{f}(x)) &= \text{local positioning cost}(\text{signal}(x_i))  \\
& \quad\quad + \text{formation cost}(\|\text{signal}(x_i) - \sbs{\sigma}{f}(x)\|), \notag
\end{align*}
where $\text{signal}(x_i)$ is a function whose analytic form is unknown. It is 
immediate to see that, in this framework, the 
requirements~\ref{requirement_1}--\ref{requirement_2} naturally apply.

\textit{(Smart grids \& Energy sharing markets)}
In a shared energy system (see Fig.~\ref{fig:applications}(a)), each agent is 
a prosumer who decides how much energy to consume or produce (inject into the 
grid). The local loss $\ftwo_i$ then depends on a local utility 
(describing, e.g., the comfort of the user), plus the monetary cost of 
electricity (which depends on the total demand), modeled by 
\( \sbs{\sigma}{f}(x) \). 
Accordingly, the cost of agent \( i \) can be modeled as:
\begin{align}\label{eq:example-energymarkets}
\ftwo_i(x_i, \sbs{\sigma}{f}(x)) = \text{local utility}(x_i) - \text{price}(\sbs{\sigma}{f}(x)) \cdot x_i.
\end{align}
In practice, however, the local utility function depends on user-specific 
parameters, such as comfort preferences (e.g., air-conditioning setpoints), 
productivity, the state of charge of local battery storage, and the degree of 
user participation in demand-response programs. 
Moreover, the energy price \( \sbs{\sigma}{f}(x) \) often follows an unknown 
dynamic pricing model determined by the grid operator, whose analytical form is 
not available to the user. In these cases, \eqref{eq:example-energymarkets} 
modifies to:
\begin{align*}
\ftwo_i(x_i, \sbs{\sigma}{f}(x)) = \text{comfort}_i(x_i) - 
\text{signal price}(\sbs{\sigma}{f}(x)) \cdot x_i,
\end{align*}
where $\text{comfort}(x_i)$ is a user-defined level of comfort (whose analytic form is unknown and user-specific), and $\text{signal price}(\sbs{\sigma}{f}(x))$ is a tariff set by the grid that the agents can observe but not model 
analytically (or differentiate). 
It is evident that \ref{requirement_1}--\ref{requirement_2} naturally capture 
this setting.

\textit{(Traffic networks)}
In a transportation network (see Fig.~\ref{fig:applications}(b)), each agent 
is a user who selects a mode of transportation (e.g., private or shared 
mobility) or a commuting route based on individual preferences. 
Depending on the collective choices of all users, the network experiences a 
certain level of congestion (e.g., delays), captured by \(\sbs{\sigma}{f}(x) \).
In turn, the transportation cost (e.g., transit time) of each user \( i \) is a 
combination of both their individual decision and the overall congestion:
\begin{align}\label{eq:example-congestion}
\ftwo_i(x_i, \sbs{\sigma}{f}(x)) = \text{travel cost}(x_i) + \text{congestion penalty}(\sbs{\sigma}{f}(x)).
\end{align}
In practice, however, congestion-related penalties are difficult to predict, 
as travel times can fluctuate rapidly, especially when many users decide to 
commute simultaneously. In other situations, where congestion penalties are 
used to model externalities (such as emissions), these costs are nearly 
impossible to anticipate: for instance, the emissions associated with a bus 
remain largely unchanged regardless of whether a single additional passenger 
boards it or not. Such costs can therefore be estimated only \emph{a posteriori}. In these cases, \eqref{eq:example-congestion} becomes:
\begin{align*} 
\ftwo_i(x_i, \sbs{\sigma}{f}(x)) &= \text{travel cost}(x_i)\\
&\quad\quad + \text{estimated congestion penalty}(\sbs{\sigma}{f}(x)).
\end{align*}
It is immediate that~\ref{requirement_1}--\ref{requirement_2} naturally 
capture this setting.

\textit{(Federated Learning with Personalized Objectives)}
In federated learning (see Fig.~\ref{fig:applications}(c)), each agent (e.g., a 
device or user) aims to collaboratively train a machine learning model 
(described by the hyperparameters \( x_i \)) based on data they would like to 
maintain private. The overall performance of the learning system depends on an 
aggregated model or statistics \( \sbs{\sigma}{f}(x) \), such as the average of 
all local models possibly evaluated on global validation data.
As a result, the objective of agent \( i \) shall balance local accuracy and 
global consistency, taking the form:
\begin{align}\label{eq:example-learning}
\ftwo_i(x_i, \sigma_f(x)) &= \text{local empirical risk}(x_i) \notag\\
&\quad\quad\quad+ \text{regularization}(\|x_i - \sbs{\sigma}{f}(x)\|).    
\end{align}
In practice, the relationship between the hyperparameters and the resulting 
validation performance is often non-analytic and noisy; 
in these cases, the requirements~\ref{requirement_1}--\ref{requirement_2} 
arise naturally as gradients are unavailable and costs are known only locally.

%%%%%%%%%%%%%%%%%%%%%%%%%%%%%%%%%%%%%%%%%%%%%%%%%%%%%%%%%%%%%%%%%%%%%%
%%%%%%%%%%%%%%%%%%%%%%%%%%%%%%%%%%%%%%%%%%%%%%%%%%%%%%%%%%%%%%%%%%%%%%%%%%%%%%%
\section{Application of the framework to solve robotic formation control problems}
\label{sec:simulations}

\begin{figure*}[t]
\centering \subfigure{\includegraphics[width=.68\columnwidth]{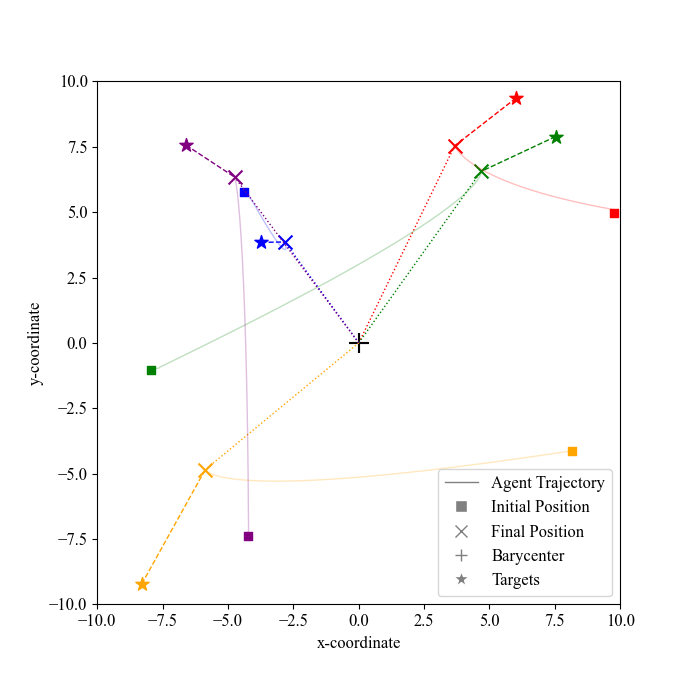}} 
\raisebox{2mm}{\centering \subfigure{\includegraphics[width=.6\columnwidth]{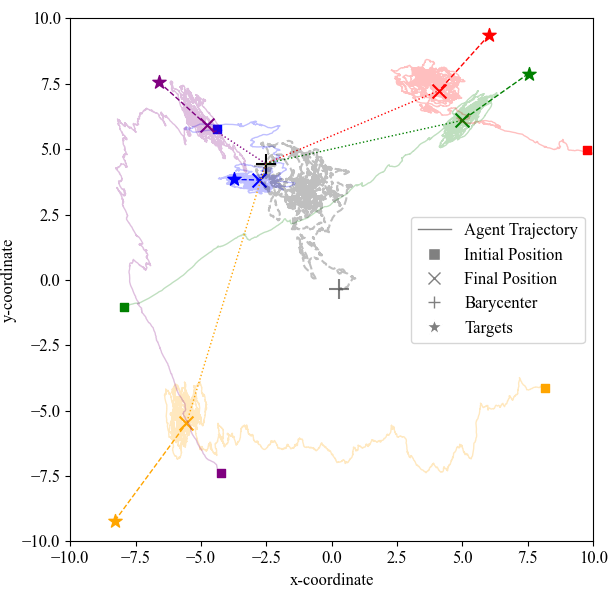}}}
\centering \subfigure{\includegraphics[width=.68\columnwidth]{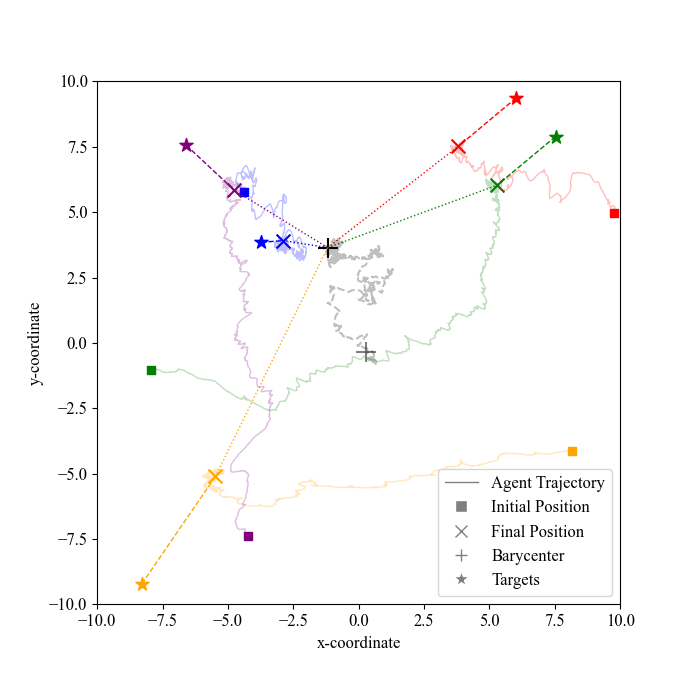}} 
\caption{
\textbf{(Left)} Robot trajectories under the algorithm of~\cite{XL-LX-YH:22}, which relies on exact gradient information. 
\textbf{(Center)} Trajectories obtained with \texttt{ARGFree} (Algorithm~\ref{alg:ARGFree}). 
\textbf{(Right)} Trajectories obtained with \texttt{ARGFree-EM} (Algorithm~\ref{alg:ARGFree_filtration}). 
Both proposed methods achieve performance comparable to the gradient-based approach while requiring no gradient information. 
Qualitatively, \texttt{ARGFree-EM} produces smoother trajectories that more closely resemble those of the gradient-based method.
}
\vspace{-.5cm}
\label{fig:sims_new}
\end{figure*}

% \begin{figure}[t]
% \centering \subfigure[]{\includegraphics[width=1\columnwidth]{paper/plots/CostComp_noNoise_final.png}} 
% % \centering \subfigure[]{\includegraphics[width=\columnwidth]{paper/plots/fig1.png}}\\
% \centering \subfigure[]{\includegraphics[width=1\columnwidth]{paper/plots/CostComp_withNoise_final.png}}\\
% \caption{Application of Algorithm~\ref{alg:ARGFree_filtration} to solve a 
% robotic formation control problem.  (a) Paths followed by the robots. (b) 
% Temporal evolution of the loss function and relative loss. All simulations are 
% averaged over Monte Carlo experiments of $10$ simulations; shaded areas 
% illustrate $\pm\sigma$ confidence intervals.
% Orange lines illustrate a comparison with the method from~\cite{XL-LX-YH:22} 
% relying on exact gradient knowledge.
% The illustration shows that our method allows the agents to reach the desired 
% configuration up to a relative error of the order $10^{-2}.$ Relative to our 
% method, approaches using exact gradients exhibit superior convergence rate and accuracy, although they require additional knowledge and cannot be implemented under requirements~\ref{requirement_1}--\ref{requirement_2}.
% % The comparison shows that using finite-different approximations in place of 
% % exact gradients corresponds to a loss of performance on both convergence rate 
% % and accuracy, but offer the advantage of requiring no gradient information 
% % while still reaching the optimizer within $10^{-2}$ relative error.
% }
% \vspace{-.5cm}
% \label{fig:sims_noiseless}
% \end{figure}

In this section, we demonstrate the applicability of the proposed methods 
through numerical simulations. 
We consider a multi-agent robotic formation control problem, 
where each agent \( i \) is a robot with position \( x_i \in \mathbb{R}^2 \) 
that aims to reach its privately known {\it target} \( r_i \in \mathbb{R}^2 \).
At the same time, the robots must maintain cohesion 
within the swarm, which is enforced by penalizing large deviations of 
\(x_i\) from the swarm's {\it barycenter} 
$\sbs{\sigma}{f}(x):= \sum_{i=1}^N x_i/N$. 
This problem can be modeled as an instance of~\eqref{eq:sigma_x} with 
$\phi_i(x_i)=x_i$ and
\begin{align*}
\ftwo_i(x_i,\sbs{\sigma}{f}(x))=\dfrac{\gamma_i}{2}\|x_i-r_i\|^2+\dfrac{1}{2}\|x_i-\sbs{\sigma}{f}(x)\|^2,
\end{align*}
where \(\gamma_i > 0\) is a weighting parameter that models the extent to which 
robot \(i\) prioritizes reaching its target position \(r_i \in \real^2\), over 
maintaining cohesion with the swarm.

For our simulations, we considered a system of \(N = 5\) agents communicating 
over an instance of a random graphs with Erd\H{o}s–R\'enyi topology, edge 
probability \(p = 0.6\), and 
uniform weights---see Assumption~\ref{as:graph}. 
We set \(\gamma_i = 2\) for all agents and generated both the desired positions 
\(r_i\) and initial positions \(x_i^0\) uniformly at random in the 
interval \([0, 10]\).
The algorithm's parameters have been chosen as follows: 
\(\alpha = 2 \cdot 10^{-3}\), \(\delta = 10^{-5}\), \(B_i\) a 
randomly generated matrix with eigenvalues in the range \((0.9, 1)\), 
\(\Sigma_u^0 = I\), and \(\Sigma_v = 0.16 I\). This configuration was selected 
after testing several alternatives, as it consistently produced consistent behavior 
and reliable performance across a variety of simulation scenarios.

\begin{figure}[t]
\centering \subfigure[]{\includegraphics[width=\columnwidth]{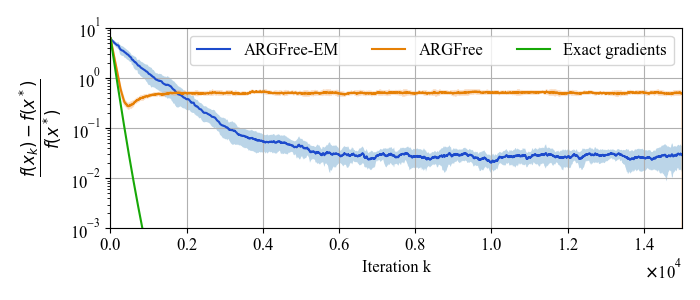}} 
\centering \subfigure[]{\includegraphics[width=\columnwidth]{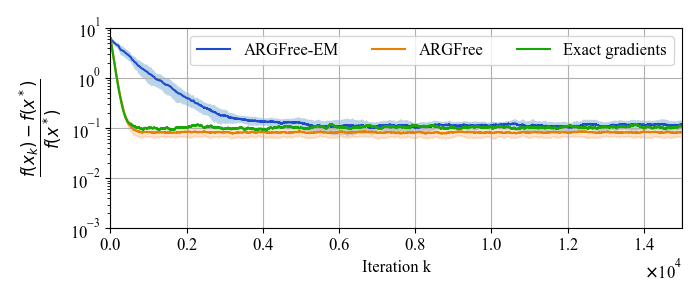}}
\caption{Evolution of the relative loss. (a) Noiseless scenario. (b) Noisy scenario. Simulations are averaged over 10 Monte Carlo runs, with shaded areas indicating $\pm\sigma$ confidence intervals. Relative to our 
method, approaches using exact gradients exhibit superior convergence rate and accuracy in the noiseless case although they require additional knowledge and cannot be implemented under requirements~\ref{requirement_1}--\ref{requirement_2}. In the presence of measurement noise, the performance of 
exact-gradient methods deteriorates significantly, whereas the proposed \texttt{ARGFree} and \texttt{ARGFree-EM} algorithms have superior accuracy.}
\label{fig:sims_loss}
\vspace{-.2cm}
\end{figure}

\subsection{Noiseless scenario}
\label{sec:simulations_noiseless}

% Simulation results from the implementation of 
% Algorithm~\ref{alg:ARGFree_filtration} are presented in 
% Fig.~\ref{fig:sims_noiseless}. 
% Fig.~\ref{fig:sims_noiseless}(a) shows the evolution of the robots' positions 
% on the plane. The trajectories demonstrate that the robots asymptotically 
% converge to configurations (indicated by `$\times$' markers) that strike a 
% tradeoff between reaching their individual target positions (shown as `$\ast$') 
% and maintaining proximity to the barycenter of the swarm (denoted by the `$+$' 
% symbol). Fig.~\ref{fig:sims_noiseless}(b)-top illustrates that the algorithm is 
% effective at minimizing the function $\fone(x)$, while 
% Fig.~\ref{fig:sims_noiseless}(b)-bottom shows that the asymptotic relative 
% error of the method is within $10^{-1}$ accuracy. 
% Fig.~\ref{fig:sims_noiseless}(b) also presents a comparison between our approach 
% (blue lines) and an exact method that harnesses exact knowledge of gradients of 
% the functions involved \(\ftwo_i(x_i, \sigma)\), taken from~\cite{XL-LX-YH:22}.
% The comparison highlights that using a finite-difference scheme with random 
% perturbations leads to a loss in performance relative to the exact method, both 
% in terms of convergence rate and asymptotic accuracy. However, this tradeoff 
% comes with the advantage of reduced informational requirements, as our method 
% does not rely on analytic gradient information---see 
% requirements~\ref{requirement_1}–\ref{requirement_2}.

Simulation results comparing the proposed algorithms against an exact gradient-based method~\cite{XL-LX-YH:22} are presented in Fig.~\ref{fig:sims_new}, Fig.~\ref{fig:sims_loss}(a), and Fig.~\ref{fig:sims_grad}(a).

Fig.~\ref{fig:sims_new} illustrates the robot trajectories. Specifically, it 
compares the exact gradient-based method from~\cite{XL-LX-YH:22} (left panel) 
with \texttt{ARGFree} (Algorithm~\ref{alg:ARGFree}) (center panel) and 
\texttt{ARGFree-EM} (Algorithm~\ref{alg:ARGFree_filtration}) (right panel). 
All three algorithms drive the robots to asymptotic configurations (marked by `$\times$') that balance reaching individual targets (`$\ast$') and maintaining proximity to the swarm barycenter (`$+$'). 
Qualitatively, \texttt{ARGFree-EM} yields smoother trajectories than \texttt{ARGFree}, more closely resembling those of the gradient-based method. As expected, \texttt{ARGFree} produces more irregular paths due to randomized perturbations used for gradient estimation; nevertheless, all methods converge to nearly identical configurations. 
Overall, this demonstrates that the proposed approach achieves performance comparable to exact gradient methods without requiring gradient information.

Fig.~\ref{fig:sims_loss}(a) shows the relative loss, confirming that our methods effectively minimize the global objective to an asymptotic accuracy of order $10^{-2}$. Fig.~\ref{fig:sims_grad}(a) corroborates this by showing the corresponding decrease in the gradient norm. The comparison highlights a clear trade-off: using a randomized finite-difference scheme incurs a penalty in convergence rate and accuracy relative to the exact method. However, this is expected and offsets the strict informational requirement of requiring analytic gradients---see 
requirements~\ref{requirement_1}–\ref{requirement_2}.

% Fig.~\ref{fig:sims_new} illustrates the paths followed by the robots. 
% Specifically, it compares the exact gradient-based method 
% from~\cite{XL-LX-YH:22} (left panel) with \texttt{ARGFree}
% (Algorithm~\ref{alg:ARGFree}) (right panel) and 
% and \texttt{ARGFree-EM} (Algorithm~\ref{alg:ARGFree_filtration}) (right panel). 
% All three algorithms drive the robots to asymptotic configurations 
% (marked by  `$\times$' symbols) that balance reaching their individual targets 
% (denoted by `$\ast$' symbols)  and maintaining proximity to the swarm 
% barycenter (denoted by the `$+$' symbol). 
% Qualitatively, \texttt{ARGFree-EM} produces smoother trajectories relative to 
% \texttt{ARGFree} that more closely resemble those of the gradient-based method.
% As expected, the trajectories generated by \texttt{ARGFree} are more irregular 
% due to the use of randomized perturbations for gradient estimation; however, 
% the robots converge asymptotically to (visually) nearly identical 
% configurations. 
% This demonstrates that the proposed approach achieves comparable performance 
% to exact gradient methods despite not requiring gradient information for its 
% implementation.

\subsection{Noisy scenario}

In practical applications, absolute positions \(x_i\) are often difficult to 
estimate reliably. Instead, agents typically rely on noisy local measurements 
(e.g., radio signal strength) at their current location; see 
Section~\ref{sec:illustrative_applications}--\textit{(Multi-agent robotics)}, for 
a more detailed discussion. The effect of noisy position measurements can be 
modeled by replacing the quantity \(x_i\) with \(w_i  x_i\) in the right-hand 
side of lines 1–7 of Algorithm~\ref{alg:ARGFree_filtration}, where 
\(w_i \sim \mathcal{N}(0, 0.2 I_2)\). These inexact values are then used to 
compute the update for robot \(i\)'s true position \(x_i^{k+1}\) in line 1 of 
the algorithm.

\begin{figure}[t]
\centering \subfigure[]{\includegraphics[width=\columnwidth]{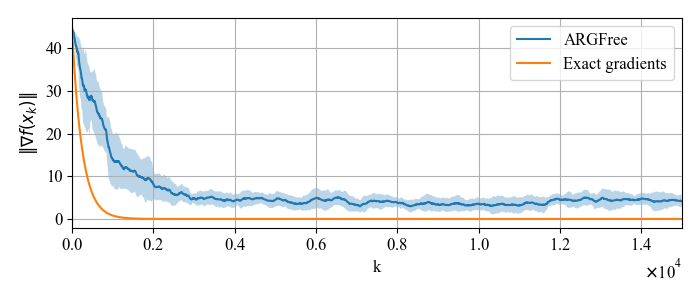}} 
\centering \subfigure[]{\includegraphics[width=\columnwidth]{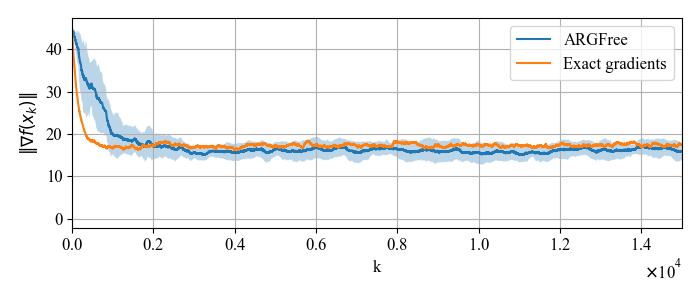}}
\caption{Evolution of the gradient norm $\|\nabla f(x_k)\|$ evaluated at the current iterates. (a) Noiseless scenario. (b) Noisy scenario. These results corroborate the relative loss findings, demonstrating the robustness of direct function evaluations against gradient-amplified noise.}
\label{fig:sims_grad}
\vspace{-.2cm}
\end{figure}

The performance under this measurement noise is shown in Fig.~\ref{fig:sims_loss}(b) for relative loss and Fig.~\ref{fig:sims_grad}(b) for the gradient norm. While the exact gradient method outperformed our approach in the noiseless setting, its performance deteriorates significantly here. This outcome can be explained by 
noting that, in the exact method, small perturbations in the estimated position 
\(x_i\) are amplified through the gradient computation. In comparison, our 
algorithm relies solely on direct evaluations of the objective function (i.e., 
measurements of position-dependent signals), which tend to be less sensitive to 
noise than their gradients. As a result, our algorithms achieve a superior asymptotic accuracy in the noisy regime, underscoring the practical robustness of the gradient-free framework.

%%%%%%%%%%%%%%%%%%%%%%%%%%%%%%%%%%%%%%%%%%%%%%%%%%%%%%%%%%%%%%%%%%%%%%
%%%%%%%%%%%%%%%%%%%%%%%%%%%%%%%%%%%%%%%%%%%%%%%%%%%%%%%%%%%%%%%%%%%%%%%%%%%%%%%
\section{Conclusions}
\label{sec:conclusion}
We have proposed \texttt{ARGFree} and \texttt{ARGFree-EM}, two distributed 
algorithms for solving aggregative cooperative optimization problems without 
requiring explicit gradient information. 
Both methods rely on randomized finite-difference approximations of the cost 
gradient. 
We established that the algorithms converge to a neighborhood of the 
optimizer, whose size can be controlled through an appropriate choice of the 
smoothing parameter. 
Compared to \texttt{ARGFree}, the \texttt{ARGFree-EM} variant achieves higher 
accuracy by incorporating an exploration momentum that mitigates rapid 
fluctuations in the exploration signal.
This work opens the opportunity for several directions of future work, 
including the use of single-point and multi-point gradient approximations, and 
adaptations of the methods for feedback optimization configurations.

\appendices

\section{Proof of Proposition~\ref{prop:tracking_variables}}
\label{sec:proof_lem_tracking_variables}

We begin by proving $\bar \sigma_k = \sbs{\sigma}{f}(x_k).$
By multiplying by $\frac{1}{N} \one^\top$ both sides of \eqref{eq:sigma} and by 
using $\frac{1}{N} \one^\top A = \frac{1}{N} \one^\top$ (from 
Assumption~\ref{as:graph}), we have:
\begin{align*}
\bar{\sigma}_{k+1}=\bar{\sigma}_k 
+ \frac{1}{N}\sum_{i=1}^N \left( \phi_i(x_i^{k+1})- \phi_i(x_i^{k}) \right).
\end{align*}
By telescoping the sum:
\begin{align*}
\bar{\sigma}_{k+1} &= \bar \sigma_0
+ \frac{1}{N}\sum_{i=1}^N \left( \phi_i(x_i^{k+1})- \phi_i(x_i^{0}) \right)
\\
&= \frac{1}{N}\sum_{i=1}^N \phi_i(x_i^{k+1}),
\end{align*}
where the last identity follows from the choice of initial conditions 
$\sigma_{i}^{0} = \phi_i(x_{i}^{0}).$ The first assertion thus follows by 
definition of $\sbs{\sigma}{f}(x)$ (see~\eqref{eq:sigma_x}).
The proof of the remaining assertions follows by iterating the argument.

\section{Proof of Lemma~\ref{lem:bound_twopoint_alg}}
\label{sec:proof_lem_bound_twopoint_alg}

We begin by recalling the following basic properties: 
(a) $\fone(x) - \fone(x^\star) \leq {\nabla \fone(x)}^T (x-x^\star)$, which follows from 
convexity of $\fone(x)$; and (b)
$\|\nabla \fone(x)\|^2 \leq 2 \Lone_1 \left(\fone(x) - \fone(x^\star)\right)$ 
which follows from Lipschitz smoothness. We have:
\begin{align*}
E_{u} & [\norm{x - \alpha g_\delta(x) - x^*}^2]\\
&= 
\norm{x - x^*}^2
-2 \alpha \left\langle E_{u}[g_\delta (x)], x-x^*\right\rangle
+\alpha^2 E_{u}[\norm{g_\delta (x)}^2]\\
&\reljust{\leq}{\eqref{eq:Nesterov35}}
\norm{x - x^*}^2
-2 \alpha \left\langle E_{u}[g_\delta (x)], x-x^*\right\rangle\\
&\quad\quad +\alpha^2 \left[ \frac{\delta^2 (n+6)^3}{2} \Lone_1^2+2(n+4)\|\nabla \fone(x)\|^2\right]\\
&\reljust{\leq}{\eqref{eq:Nesterov21},(a)}
\norm{x - x^*}^2
-2 \alpha \left( \fone_\delta(x) - \fone_\delta(x^\star) \right)\\
&\quad\quad +\alpha^2 \left[ \frac{\delta^2 (n+6)^3}{2} \Lone_1^2+2(n+4)\|\nabla \fone(x)\|^2\right]\\
&\reljust{\leq}{\eqref{eq:Nesterov11},(b)}
\norm{x - x^*}^2
-2 \alpha (\fone(x) - \fone_\delta(x^*))\\ 
&\quad\quad +\alpha^2 \left[ \frac{\delta^2 (n+6)^3}{2} \Lone_1^2+4(n+4)\Lone_1(\fone(x) - \fone(x^*))\right]\\
&\reljust{\leq}{\eqref{eq:Nesterov19}}
\norm{x - x^*}^2
-2 \alpha \left( 1- 2\alpha(n+4)\Lone_1 \right)(\fone(x) - \fone(x^*))\\ 
&\quad\quad + \delta^2 n \alpha \Lone_1 + \frac{\delta^2 (n+6)^3}{2} \alpha^2 \Lone_1^2\\
&\leq 
\norm{x - x^*}^2
- \alpha \mu (1-2\alpha(n+4)\Lone_1)\norm{x  - x^*}^2\\ 
&\quad\quad + \delta^2 n \alpha \Lone_1 + \frac{\delta^2 (n+6)^3}{2} \alpha^2 \Lone_1^2,
\end{align*}
where for the last inequality we used
$\fone(x)-\fone\left(x^*\right) \geq \mu/2 \left\|x-x^*\right\|^2.$
This establishes that
\begin{align}
\label{eq:contraction_squared}
E_{u} [\norm{x - \alpha g_\delta(x)- x^*}^2] \leq 
(1-\beta_1^\alpha) \norm{x - x^*}^2 + (\beta_2^\alpha)^2.
\end{align}
To conclude, notice that
\begin{align*}
E_{u} [\norm{x - \alpha g_\delta(x)- x^*}] &\leq    
\left(E_{u} [\norm{x - \alpha g_\delta(x)- x^*}^2]\right)^{1/2}\\
&\reljust{\leq}{\eqref{eq:contraction_squared}}
\sqrt{1-\beta_1^\alpha} \norm{x - x^*} + \beta_2^\alpha,
\end{align*}
where the first bound follows by Jensen’s inequality, and the second one from $\sqrt{a+b}\leq \sqrt{a}+\sqrt{b}$ for $a,b \geq 0.$

\bibliographystyle{myIEEEtran}
\bibliography{BIB/alias,BIB/full_GB,BIB/GB,BIB/AM_bib}

\end{document}